\documentclass[12pt,english]{article}
\usepackage[T1]{fontenc}
\usepackage[latin9]{inputenc}
\usepackage[letterpaper]{geometry}
\usepackage{babel}
\usepackage{prettyref}
\usepackage{amsthm}
\usepackage{amsmath}
\usepackage{amssymb}
\usepackage{graphicx}
\usepackage[all]{xy}
\usepackage{epstopdf}
\usepackage[unicode=true,pdfusetitle,
 bookmarks=true,bookmarksnumbered=false,bookmarksopen=false,
 breaklinks=false,pdfborder={0 0 1},backref=section,colorlinks=false]
 {hyperref}

\makeatletter

\providecommand{\tabularnewline}{\\}

\numberwithin{equation}{section}
\numberwithin{figure}{section}
\theoremstyle{plain}
\newtheorem{thm}{\protect\theoremname}[section]
  \theoremstyle{remark}
  \newtheorem{rem}[thm]{\protect\remarkname}
  \theoremstyle{definition}
  \newtheorem{defn}[thm]{\protect\definitionname}
  \theoremstyle{plain}
  \newtheorem{conjecture}[thm]{\protect\conjecturename}
  \theoremstyle{plain}
  \newtheorem{cor}[thm]{\protect\corollaryname}
  \theoremstyle{plain}
  \newtheorem{prop}[thm]{\protect\propositionname}
  \theoremstyle{remark}
  \newtheorem{claim}[thm]{\protect\claimname}
  \theoremstyle{plain}
  \newtheorem{lem}[thm]{\protect\lemmaname}

\usepackage[leftcaption]{sidecap}

\sidecaptionvpos{figure}{c}

\newcommand{\FigBesBeg}[1][1.0]{%
 \let\MyFigure\figure
 \let\MyEndfigure\endfigure
 \renewenvironment{figure}[1]{\begin{SCfigure}[#1]##1}{\end{SCfigure}}}

\newcommand{\FigBesEnd}{%
 \let\figure\MyFigure
 \let\endfigure\MyEndfigure}

\usepackage{amsfonts}
\usepackage{amsthm}
\usepackage[all]{xy}
\usepackage{stackrel}
\usepackage{appendix}

\newcommand{\f}{\phi}
\newcommand{\F}{\mathbf{F}}
\newcommand{\hF}{{\hat{\F}}}

\renewcommand{\G}{\Gamma}
\newcommand{\hG}{{\hat{\Gamma}}}
\renewcommand{\U}{\Upsilon}
\newcommand{\ff}{\stackrel{*}{\le}}
\newcommand{\fg}{\le_{fg}}
\newcommand{\covers}{\stackrel{X}{\twoheadrightarrow}}
\newcommand{\D}{\mathcal{D}}
\newcommand{\la}{\langle}
\newcommand{\ra}{\rangle}

\renewcommand{\O}{\mathcal{O}}
\renewcommand{\a}{\alpha}
\newcommand{\Co}{{Comp}}

\theoremstyle{definition}
\newtheorem{rem}[thm]{\protect\remarkname}

\theoremstyle{definition}
\newtheorem{defn}[thm]{\protect\definitionname}

\theoremstyle{plain}
\newtheorem{claim}[thm]{\protect\claimname}

\makeatother

  \providecommand{\claimname}{Claim}
  \providecommand{\conjecturename}{Conjecture}
  \providecommand{\corollaryname}{Corollary}
  \providecommand{\definitionname}{Definition}
  \providecommand{\lemmaname}{Lemma}
  \providecommand{\propositionname}{Proposition}
  \providecommand{\remarkname}{Remark}
\providecommand{\theoremname}{Theorem}

\begin{document}

\title{Primitive Words, Free Factors and Measure Preservation}

\author{Doron Puder%
\thanks{Supported by Advanced ERC Grant 247034 of Aner Shalev, and by Adams
Fellowship Program of the Israel Academy of Sciences and Humanities.%
}\\
Einstein Institute of Mathematics \\
 Hebrew University, Jerusalem\\
 \texttt{doronpuder@gmail.com}}
\maketitle
\begin{abstract}
Let $\F_{k}$ be the free group on $k$ generators. A word $w\in\F_{k}$
is called primitive if it belongs to some basis of $\F_{k}$. We investigate
two criteria for primitivity, and consider more generally, subgroups
of $\F_{k}$ which are free factors.

The first criterion is graph-theoretic and uses Stallings core graphs:
given subgroups of finite rank $H\le J\le\F_{k}$ we present a simple
procedure to determine whether $H$ is a free factor of $J$. This
yields, in particular, a procedure to determine whether a given element
in $\F_{k}$ is primitive.

Again let $w\in\F_{k}$ and consider the word map $w:G\times\ldots\times G\to G$
(from the direct product of $k$ copies of $G$ to $G$), where $G$
is an arbitrary finite group. We call $w$ \emph{measure preserving}
if given uniform measure on $G\times\ldots\times G$, $w$ induces
uniform measure on $G$ (for every finite $G$). This is the second
criterion we investigate: it is not hard to see that primitivity implies
measure preservation and it was conjectured that the two properties
are equivalent. Our combinatorial approach to primitivity allows us
to make progress on this problem and in particular prove the conjecture
for $k=2$.

It was asked whether the primitive elements of $\F_{k}$ form a closed
set in the profinite topology of free groups. Our results provide
a positive answer for $\F_{2}$.

\end{abstract}

\section{Introduction}

An element $w$ of a free group $J$ is called \emph{primitive} if
it belongs to some basis (free generating set) of $J$. When $J$
is given with a basis $X$, this is equivalent to the existence of
an automorphism of $J$ which sends $w$ to a given element of $X$.

The notion of primitivity has a natural extension to subgroups in
the form of free factors. Let $H$ be a subgroup of the free group
$J$ (in particular, $H$ is free as well). We say that $H$ is a
\emph{free factor} of $J$ and denote $H\ff J$, if there is another
subgroup $H'\le J$ such that $H*H'=J$. Equivalently, $H\ff J$ if
every basis of $H$ can be extended to a basis of $J$. (This in turn
is easily seen to be equivalent to the condition that {\em some}
basis of $H$ extends to a basis of $J$).

Let $\F_{k}$ be the free group on $k$ generators with a fixed basis
$X=\{x_{1},\ldots,x_{k}\}$. We study finitely generated subgroups
of $\F_{k}$ (denoted $H\fg\F_{k}$) and relations among them using
core graphs, also known as Stallings' graphs (See~\cite{Sta83}.
Actually our definition is a bit different than Stalling's, see below).
Associated with every subgroup $H\le\F_{k}$ is a directed, pointed,
edge-labeled graph denoted $\G_{X}(H)$. Edges are labeled by the
elements of the given basis $X=\{x_{1},\ldots,x_{k}\}$ of $\F_{k}$.
A full definition appears in Section \ref{sec:core-graphs}, but we
illustrate the concept in Figure \ref{fig:first_core_graph}. It shows
the core-graph of the subgroup of $\F_{2}$ generated by $x_{1}x_{2}x_{1}^{\;-1}x_{2}^{\;-1}$
and $x_{2}x_{1}^{\;2}$.

\FigBesBeg 
\begin{figure}[h]
\centering{}%
\begin{minipage}[t]{0.5\columnwidth}%
\[
\xymatrix{\otimes\ar[rr]^{1} &  & \bullet\\
\\
\bullet\ar[rr]^{1}\ar[uu]_{2} &  & \bullet\ar[uull]_{1}\ar[uu]_{2}
}
\]
\end{minipage}\caption{\label{fig:first_core_graph} The core graph $\G_{X}\left(H\right)$
where $H=\left\langle x_{1}x_{2}^{-1}x_{1},x_{1}^{-2}x_{2}\right\rangle \leq\F_{2}$.}
\end{figure}
\FigBesEnd 

Core graphs are a key tool in the research of free groups, and are
both used for proving new results and for introducing simple proofs
to known results (see, for instance, \cite{KM02,MVW07}, for a survey
of many such results and for further references).

A central new ingredient of our work is a new perspective on core
graphs. There is a naturally defined notion of quotient on such graphs
(see Section \ref{sec:dag}). In particular, we introduce in Section
\ref{sec:dag} the notion of \emph{immediate quotients}. This in turn
yields a directed graph whose vertices are all core graphs of finitely
generated subgroups of $\F_{k}$ (w.r.t. the fixed basis $X$). A
directed edge in this graph stands for the relation of an immediate
quotient. This is a directed acyclic graph (DAG) i.e., it contains
no directed cycles. As always, reachability in a DAG induces a distance
function between vertices. Namely $\rho_{X}(x,y)$ is the shortest
length of a directed path from $x$ to $y$. We mention that the transitive
closure of the immediate quotient relation is the relation {}``being
a quotient of'' which is a partial order (a lattice, in fact) on
all core graphs of f.g. subgroups of $\F_{k}$. The following theorem
gives a simple criterion for free factorness in terms of this distance:
\begin{thm}
\label{thr:rho=00003Drk-rk_iff_ff} Let $H,J\fg\F_{k}$, and assume
$\G_{X}(J)$ is a quotient of $\G_{X}(H)$. Then $H\ff J$ if and
only if 
\[
\rho_{X}(H,J)=rk(J)-rk(H)
\]
 
\end{thm}
We note that $\rho_{X}(\cdot,\cdot)$ can be explicitly computed,
and this theorem thus yields automatically an algorithm to determine,
for two given $H,J\fg\F_{k}$ whether $H$ is a free factor of $J$.
In particular, it can serve to detect primitive words (see Appendix
\ref{sec:algo}). More generally, for any f.g.~free groups $H\leq J$,
this theorem can serve to detect the minimal number of complementary
generators needed to obtain $J$ from $H$ (Corollary \prettyref{cor:min_num_of_complemetary}).

In fact, the core graph of every $H\fg\F_{k}$ has finitely many quotients
(or reachable vertices). This set is also known in the literature
as the \emph{fringe} of $H$ (see, e.g. \cite{MVW07}). For example,
Figure \ref{fig:lattice} shows the fringe of the subgroup $H=\la[x_{1},x_{2}]\ra$.
The difference in ranks between $H$ and $\F_{2}$ is 1. However,
the distance between the corresponding core graphs in the fringe is
2. This proves that $H$ is not a free factor of $\F_{2}$, or equivalently
that $[x_{1},x_{2}]$ is not primitive. We elaborate more in Appendix
\ref{sbs:examples}. 
\begin{rem}
We stress that there are other graph-theoretic algorithms to detect
free factors and primitive words, including simplifications of the
seminal Whitehead algorithm (the algorithm first appeared in \cite{Whi36a,Whi36b},
for its graph-theoretic simplifications see \cite{Ger84,Sta99}).
Our approach, however, is very different and does not rely on Whitehead
automorphisms. We elaborate more on this in Appendix \ref{sec:algo}.
\end{rem}
Theorem \ref{thr:rho=00003Drk-rk_iff_ff} is also used for the other
concept we study here, that of measure preservation of word maps.
Associated with every $w\in\mathbf{F}_{k}$ is a \emph{word map}.
We view $w$ as a word in the letters of the basis $X$. For every
group $G$, this mapping which we also denote by $w$ maps $\underbrace{G\times G\times\cdots\times G}_{k}\longrightarrow G$
as follows: It maps the $k$-tuple $(g_{1},\ldots,g_{k})$ to the
element $w(g_{1},\ldots,g_{k})\in G$, where $w(g_{1},\ldots,g_{k})$
is the element obtained by replacing $x_{1},\ldots,x_{k}$ with $g_{1},\ldots,g_{k}$
(respectively) in the expression for $w$, and then evaluating this
expression as a group element in $G$.

During the last years there has been a great interest in word maps
in groups, and extensive research was conducted (see, for instance,
\cite{Sha09}, \cite{LSh09}; for a recent book on the topic see \cite{Seg09}).
Our focus here is on the property of measure preservation: We say
that the word $w$ preserves measure with respect to a finite group
$G$ if when $k$-tuples of elements from $G$ are sampled uniformly,
the image of the word map $w$ induces the uniform distribution on
$G$. (In other words, all fibers of the word map have the same size).
We say that $w$ is \emph{measure preserving} if it preserves measure
with respect to \emph{every} finite group $G$.

This concept was investigated in several recent works. See for example
\cite{LSh08} and \cite{GSh09}, where certain word maps are shown
to be almost measure preserving, in the sense that the distribution
induced by $w$ on finite simple groups $G$ tends to uniform, say,
in $L_{1}$ distance, when $|G|\to\infty$.

Measure preservation can be equivalently defined as follows: fix some
finite group $G$, and select a homomorphism $\alpha_{G}\in Hom(\F_{k},G)$
uniformly at random. A homomorphism from a free group is uniquely
determined by choosing the images of the elements of a basis, so that
every homomorphism is chosen with probability $1/|G|^{k}$. We then
say that $w\in\F_{k}$ is measure preserving if for every finite group
$G$ and a random homomorphism $\alpha_{G}$ as above, $\alpha_{G}(w)$
is uniformly distributed over $G$.

We note that there is a stronger condition of measure preservation
on a word $w$ that is discussed in the literature. In this stronger
condition we consider the image of $w$ over the broader class of
compact groups $G$ w.r.t. their Haar measure. Our results make use
only of the weaker condition that involves only finite groups.

Measure preservation can also be defined for f.g. subgroups. 
\begin{defn}
For $H\fg\F_{k}$ we say that $H$ is \emph{measure preserving} iff
for any finite group $G$ and $\alpha_{G}\in Hom(\F_{k},G)$ a randomly
chosen homomorphism as before, $\alpha_{G}|_{H}$ is uniformly distributed
in $Hom(H,G)$. 
\end{defn}
In particular, $1\ne w\in\F_{k}$ is measure preserving iff $\la w\ra$
is measure preserving.\\

It is easily seen that primitivity or free factorness yield measure
preservation. The reason is that as mentioned, a homomorphism in $Hom(\F_{k},G)$
is completely determined by the images of the elements of a basis
of $\F_{k}$, which can be chosen completely arbitrarily and independently.

Several authors have conjectured that the converse is also true:
\begin{conjecture}
\label{conj:prim<=00003D>m.p} For every $w\in\F_{k}$, 
\[
w\textrm{ is primitive }\Longleftrightarrow w\textrm{ is measure preserving }
\]
 More generally, for $H\fg\F_{k}$, 
\[
~~~~~~~~~H\ff\F_{k}\Longleftrightarrow H\textrm{ is measure preserving }
\]
 
\end{conjecture}
From private conversations we know that this has occurred to the following
mathematicians and discussed among themselves: T. Gelander, A. Shalev,
M. Larsen and A. Lubotzky. The question was mentioned several times
in the Einstein Institute Algebra Seminar. This conjecture was independently
raised in~\cite{LP10}%
\footnote{It is interesting to note that there is an easy abelian parallel to
Conjecture \ref{conj:prim<=00003D>m.p}: A word $w\in\F_{k}$ is primitive,
i.e. belongs to a basis, in $\mathbb{Z}^{k}\cong\F_{k}/\F'_{k}$ iff
for any group $G$ the associated word map is surjective. See \cite{Seg09},
Lemma 3.1.1.%
}.

Here we prove a partial result:
\begin{thm}
\label{ther:m.p.=00003D=00003D>prim._when_r>=00003Dk-1} Let $H\fg\F_{k}$
have rank $\ge k-1$. Then, 
\[
~~~~~~~~~H\ff\F_{k}\Longleftrightarrow H\textrm{ is measure preserving }
\]
 In particular, for every $w\in\F_{2}$: 
\[
w\textrm{ is primitive }\Longleftrightarrow w\textrm{ is measure preserving }
\]
 
\end{thm}
The proof of this result relies, inter alia, on Theorem \ref{thr:rho=00003Drk-rk_iff_ff}.
Note that a set of $k-1$ elements $w_{1},\ldots,w_{k-1}\in\F_{k}$
can be extended to a basis iff it is a free set that generates a free
factor. Thus, the result for subgroups can also be stated for finite
subsets as follows: Let $r\ge k-1$. A set $\{w_{1},\ldots,w_{r}\}\subset\F_{k}$
can be extended to a basis iff for every finite group $G$ and random
homomorphism $\alpha_{G}$ as above, the $r$-tuple $\left(\alpha_{G}(w_{1}),\ldots,\alpha_{G}(w_{r})\right)$
is uniformly distributed in $G^{r}$, the direct product of $r$ copies
of $G$.\\

There is an interesting connection between this circle of ideas and
the study of profinite groups. For example, an immediate corollary
of Theorem \ref{ther:m.p.=00003D=00003D>prim._when_r>=00003Dk-1}
is that 
\begin{cor}
\label{cor:prim_is_closed} The set of primitive elements in $\F_{2}$
is closed in the profinite topology. 
\end{cor}
We discuss this corollary and other related results in Section \ref{sec:profinite}.\\

In order to prove Conjecture \ref{conj:prim<=00003D>m.p}, one needs
to find for every non-primitive word $w\in\F_{k}$, some witness finite
group $G$ with respect to which $w$ is not measure preserving. Our
witnesses are always the symmetric groups $S_{n}$.

It is conceivable that our method of proof for Theorem \ref{ther:m.p.=00003D=00003D>prim._when_r>=00003Dk-1}
is powerful enough to establish Conjecture \ref{conj:prim<=00003D>m.p}.
We define two categorizations of elements (and of f.g. subgroups)
of free groups $\pi(\cdot)$ and $\phi(\cdot)$. They map every free
word and free subgroup into $\{0,1,2,3,\ldots\}\cup\{\infty\}$. We
believe these two maps are in fact identical. This, if true, yields
the general conjecture. Presently we can show that they are equivalent
under certain conditions, and this yields our partial result.

The first categorization is called \emph{the primitivity rank}. It
is a simple fact that if $w\in\F_{k}$ is primitive, then it is also
primitive in every subgroup of $\F_{k}$ containing it (see Claim
\ref{clm:free_factors}). However, if $w$ is not primitive in $\F_{k}$,
it may be either primitive or non-primitive in subgroups containing
it. But what is the smallest rank of a subgroup in which we can realize
$w$ is not primitive? Informally, how far does one have to search
in order to establish that $w$ is \emph{not} primitive in $\F_{k}$?
Concretely:
\begin{defn}
\label{def:primitivity-rank} The \textbf{primitivity rank} of $w\in\F_{k}$,
denoted $\pi(w)$, is 
\[
\pi(w)=min\left\{ rk(J)~\Big|~\begin{gathered}w\in J\le\F_{k}~s.t.\\
w\textrm{ is \textbf{not} primitive in \ensuremath{J}.}
\end{gathered}
\right\} 
\]
 If no such $J$ exists, $\pi(w)=\infty$. A subgroup $J$ for which
the minimum is obtained is called \textbf{$w$-critical}. 
\end{defn}
This extends naturally to subgroups. Namely, 
\begin{defn}
\label{def:primitivity-rank-subgp} For $H\fg\F_{k}$, the primitivity
rank of $H$ is 
\[
\pi(H)=min\left\{ rk(J)~\Big|~\begin{gathered}H\le J\le\F_{k}~s.t.\\
H\textrm{ is \textbf{not} a free factor of \ensuremath{J}.}
\end{gathered}
\right\} 
\]
 Again, if no such $J$ exists, $\pi(H)=\infty$. A subgroup $J$
for which the minimum is obtained is called \textbf{$H$-critical}. 
\end{defn}
For instance, $\pi(w)=1$ if and only if $w$ is a proper power of
another word (i.e. $w=v^{d}$ for some $v\in\mathbf{F}_{k}$ and $d\ge2$).
In Section \ref{sec:prim-rank} we show (Corollary \ref{cor:image-of-pi})
that in $\F_{k}$ the primitivity rank takes values only in $\{0,1,2,\ldots,k\}\cup\{\infty\}$
(the only word $w$ with $\pi(w)=0$ is $w=1$). Lemma \ref{lem:primitive-gets-infty}
shows, moreover, that $\pi(w)=\infty$ ($\pi(H)=\infty$, resp.) iff
$w$ is primitive ($H\ff\F_{k}$). Finally Lemma \ref{lem:additivity}
yields that $\pi$ can take on every value in $\{0,\ldots,k\}$. For
example, if $\F_{k}$ is given with some basis $X=\{x_{1},\ldots,x_{k}\}$
then for every $1\le d\le k$, $\pi(x_{1}^{\;2}\ldots x_{d}^{\;2})=d$.
It is interesting to mention that $\pi(H)$ also generalizes the notion
of \emph{compressed} subgroups, as appears, e.g., in \cite{MVW07}:
a subgroup $H\fg\F_{k}$ is compressed iff $\pi(H)\geq rk(H)$. \\

The second categorization of sets of formal words has its roots in
\cite{Nic94} and more explicitly in \cite{LP10}. It concerns homomorphisms
from $\F_{k}$ to the symmetric groups $S_{n}$, and more concretely
the probability that $1$ is a fixed point of the permutation $w(\sigma_{1},\ldots,\sigma_{k})$
for some $w\in\F_{k}$ when $\sigma_{1},\ldots,\sigma_{k}\in S_{n}$
are chosen randomly with uniform distribution. More generally, for
a subgroup $H\fg\F_{k}$ we study the probability that $1$ is a common
fixed point of (the permutations corresponding to) all elements in
$H$. We ask how much this probability deviates from the corresponding
probability in the case of measure preserving subgroups, i.e. from
$\frac{1}{n^{rk(H)}}$. (We continue the presentation for subgroups
only. This clearly generalizes the case of a word: for every word
$w\ne1$ consider the subgroup $\la w\ra$.)

Formally, for $H\fg\F_{k}$ we define the following function whose
domain is all integers $n\ge1$ where $\alpha_{n}\in Hom(\F_{k},S_{n})$
is a random homomorphism with uniform distribution: 
\begin{equation}
\Phi_{H}(n)=Prob\big[\forall w\in H~~~\alpha_{n}(w)(1)=1\big]-\frac{1}{n^{rk(H)}}\label{eq:Phi_def}
\end{equation}
 \\
 Clearly, if $H$ is measure preserving, then $\Phi_{H}$ vanishes
for every $n\ge1$.

Nica~\cite{Nic94} showed that for a fixed word $w\ne1$ and large
enough $n$, it is possible to express $\Phi_{w}(n)$ (=$\Phi_{\la w\ra}(n)$)
as a rational function in $n$. We show below that this is easily
extended to apply to $\Phi_{H}(n)$ for arbitrary $H\fg\F_{k}$. Nica's
clever observation was used in \cite{LP10} to introduce a new categorization
of free words, denoted $\f(\cdot)$, which, like $\pi(\cdot)$, associates
a non-negative integer or $\infty$ to every formal word (note that
in \cite{LP10} the notion of primitive words has a different meaning
than in the current paper). This categorization can also be extended
to arbitrary finitely generated subgroups of $\F_{k}$. More specifically,
it is shown in Section \ref{sec:phi} that for every $H\fg\F_{k}$
and $n$ large enough (say, at least the number of vertices in the
core graph of $H$), we have 
\begin{equation}
\Phi_{H}(n)=\sum_{i=0}^{\infty}a_{i}(H)\frac{1}{n^{i}}\label{eq:Phi_series}
\end{equation}
 where the coefficients $a_{i}(H)$ are integers depending only on
$H$. We define $\phi(H)$ as follows: 
\begin{equation}
\phi(H):=\left\{ \begin{array}{ll}
\textrm{the smallest integer \ensuremath{i}with }a_{i}(H)\neq0 & \textrm{ if }\Phi_{H}(n)\not\equiv0\\
\infty & \textrm{ if }\Phi_{H}(n)\equiv0
\end{array}\right.\label{eq:phi}
\end{equation}
 \\

Thus, $\phi(H)$ measures to what extent the probability that $1$
is a common fixed point of $H$ differs from $\frac{1}{n^{rk(H)}}$,
the corresponding probability if $H$ were measure preserving. The
higher $\phi(H)$ is, the closer the probability is asymptotically
to $\frac{1}{n^{rk(H)}}$. If $H$ is a measure preserving subgroup,
then $\phi(H)=\infty$. \\

As it turns out there is a strong connection between $\pi(H)$ and
$\phi(H)$. Already Nica's result can be interpreted in the language
of $\phi(\cdot)$ to say that $\phi(w)=1$ iff $w$ is a power, that
is iff $\pi(w)=1$. But the connection goes deeper. In proving this,
we calculate these functions using the core graph of $H$ and its
quotients. It turns out that both $\pi(H)$ and $\f(H)$ can be computed
explicitly via the subgraph of the DAG induced by all descendants
of $\G_{X}(H)$.

In the calculation of $\phi(H)$ we use the core graph $\G_{X}(H)$
and its quotients to partition the event that $1$ is a common fixed
point of $\alpha_{n}(w)$ of each $w\in H$ (see Section \ref{sec:phi}).
\\

Fortunately, the same core graph and quotients can also be used to
find the primitivity rank $\pi(H)$, as shown in Section \ref{sec:prim-rank}.
Lemma \ref{lem:pi=00003Dpi'} shows that all $H$-critical subgroups
(see Definition \ref{def:primitivity-rank-subgp}) are always represented
in the fringe (set of quotients) of $H$. Theorem \ref{thr:rho=00003Drk-rk_iff_ff}
then shows directly how to calculate $\pi(H)$ using the fringe. \\

We show that under certain conditions, the two categorizations $\pi(\cdot)$
and $\phi(\cdot)$ indeed coincide.
\begin{prop}
\label{prop:phi=00003Dpi} Let $H\fg\F_{k}$. Then for every $i\le rk(H)+1$, \end{prop}
\begin{enumerate}
\item $\pi(H)=i~~\Longleftrightarrow~~\f(H)=i$ 
\item Moreover, if $\pi(H)=\phi(H)=i$ then $a_{i}(H)$ equals the number
of $H$-critical subgroups of $\F_{k}$. 
\end{enumerate}
The second part of this proposition is in fact a generalization of
a result of Nica. For a single element $w\in\F_{k}$ which is a proper
power, namely $\pi(w)=\phi(w)=1$, let $w=u^{d}$ with $d$ maximal
(so $u$ is not a proper power). Let $M$ denote the number of divisors
of $d$. It is not hard to see that the number of $w$-critical subgroups
of $\F_{k}$ equals $M-1$: these subgroups are exactly $\la u^{m}\ra$
for every $1\le m<d$ such that $m|d$. This shows that the average
number of fixed points in the permutation $\alpha_{n}(w)$ goes to
$M$ as $n\to\infty$. This corresponds to Corollary $1.3$ in \cite{Nic94}
(for the case $L=1$)%
\footnote{Nica's result was more general in a different manner: it involved
the distribution of the number of $L$-cycles in the random permutation
$\alpha_{n}(w)$, for any fixed $L$. He showed that as $n\to\infty$,
the limit distribution depends only on $d$, where $w=u^{d}$ as above.%
}.

The connection between $\pi(\cdot)$ and $\phi(\cdot)$ goes beyond
the cases stated in Proposition \ref{prop:phi=00003Dpi}. To start
off, if $\pi(H)=\infty$, then $H\ff\F_{k}$ and therefore $H$ is
measure preserving, and thus $\phi(H)=\infty$. In addition, Lemma
\ref{lem:additivity} states that both $\pi(\cdot)$ and $\f(\cdot)$
are additive with respect to concatenation of words on disjoint letter
sets. Namely, if the words $w_{1},w_{2}\in\F_{k}$ have no letters
in common then $\pi(w_{1}w_{2})=\pi(w_{1})+\pi(w_{2})$ and $\phi(w_{1}w_{2})=\phi(w_{1})+\phi(w_{2})$.
Moreover, if the disjoint $w_{1}$ and $w_{2}$ satisfy both parts
of Proposition \ref{prop:phi=00003Dpi} then so does their concatenation
$w_{1}w_{2}$.

In view of this discussion, the following conjecture suggests itself
quite naturally:
\begin{conjecture}
\label{conj:phi=00003Dpi} \end{conjecture}
\begin{enumerate}
\item For every $H\fg\F_{k}$ 
\[
\pi(H)~~=~~\f(H)
\]
 
\item Moreover, $a_{\phi(H)}(H)$ equals the number of $H$-critical subgroups
of $\F_{k}$. 
\end{enumerate}
Specifically, for a single word $w$, Proposition \ref{prop:phi=00003Dpi}
states that for $i=0,1,2$, $\pi(w)=i\Leftrightarrow\f(w)=i$. As
mentioned, the possible values of $\pi(H)$ are $\{0,1,2,\ldots,k\}\cup\{\infty\}$,
and $\pi(H)=\infty$ iff $H\ff\F_{k}$. We also have $\pi(H)=\infty\Rightarrow\f(H)=\infty$
(a free factor subgroup is measure preserving). Thus, when $rk(H)\ge k-1$,
the value of $\pi(H)$ uniquely determines $\f(H)$ and the two values
coincide. In other words, when $rk(H)\ge k-1$ 
\[
\pi(H)=\phi(H).
\]

This shows, in turn, that when $H$ is measure preserving, we have
$\pi(H)=\phi(H)=\infty$, and so $H$ is a free factor. This yields
Theorem \ref{ther:m.p.=00003D=00003D>prim._when_r>=00003Dk-1}. The
same argument shows that Conjecture \ref{conj:prim<=00003D>m.p} follows
from part $(1)$ of Conjecture \ref{conj:phi=00003Dpi} and suggests,
in particular, a general strategy towards proving Conjecture \ref{conj:prim<=00003D>m.p}.\\

As an aside, the second parts of Proposition \ref{prop:phi=00003Dpi}
and Conjecture \ref{conj:phi=00003Dpi} say something interesting
on the average number of fixed points in the random permutation $\alpha_{n}(w)$.
We conjecture that for every $w$ and for large enough $n$, this
average is at least $1$. In other words, among the family of distributions
of $S_{n}$ induced by free words, a random uniformly chosen permutation
has the least average number of fixed points. This point is further
elaborated in Section \ref{sec:fixed_points}. \\

At this point we should clarify the relation of these results and
some of what we did in \cite{LP10}. There we introduced $\beta(\cdot)$
- yet another categorization of formal words. Just like $\phi(\cdot)$
and $\pi(\cdot)$ it maps every formal word to a non-negative integer
or $\infty$. As it turns out, $\pi(\cdot)$ and $\beta(\cdot)$ coincide.
This follows from Theorem \ref{thr:rho=00003Drk-rk_iff_ff} and from
Section \ref{sec:prim-rank}. The definition of $\pi(\cdot)$ is simpler
and more elegant than the original definition of $\beta(\cdot)$.
As shown in \cite{LP10} for $i=0,1$, $\phi(w)=i\iff\beta(w)=i$.
A partial proof was given there as well for the case $i=2$. In Section
\ref{sec:proofs} we complete the argument for $i=2$ and generalize
it to prove Proposition \ref{prop:phi=00003Dpi}. \\

The paper is arranged as follows. In section \ref{sec:core-graphs}
we introduce the notions of core graphs, their morphisms and their
quotients. In Section \ref{sec:dag} we introduce our new perspective
on core graphs, including the notion of immediate quotients and the
mentioned DAG, and then prove Theorem \ref{thr:rho=00003Drk-rk_iff_ff}.
In Section \ref{sec:prim-rank} we analyze the primitivity rank of
any $H\fg\F_{k}$ and show how it can be computed from the quotients
of $\G_{X}(H)$ in the DAG of finite rank subgroups of $\F_{k}$.
Section \ref{sec:phi} is devoted to proving that $\phi(H)$ is well
defined and can be indeed computed from the same descendants of $\G_{X}(H)$.
In Section \ref{sec:proofs} we establish the results connecting $\phi(\cdot)$
and $\pi(\cdot)$, culminating in the proof of Theorem \ref{ther:m.p.=00003D=00003D>prim._when_r>=00003Dk-1}.
The concluding sections are devoted to two different consequences
of the main results: the characterization of elements of $\F_{k}$
which are primitive in its profinite completion (Section \ref{sec:profinite})
and the possible values of the average number of fixed points in the
image of a word map on $S_{n}$ (Section \ref{sec:fixed_points}).
The discussion in the three appendices is not necessary for the main
results of this paper, but it does, in our view, complete the picture.
In particular, we illustrate in Appendix \ref{sec:algo} the algorithm
to detect free factor subgroups.

\section{Core Graphs and their Quotients}

\label{sec:core-graphs} All groups that appear here are subgroups
of $\F_{k}$, the free group with a given basis $X=\{x_{1},\ldots,x_{k}\}$.
Some of the relations we consider depend on the choice of the basis.
We first describe core-graphs, which play a crucial role in this paper.

\subsection{Core Graphs}

\label{sbs:core-graphs} Associated with every subgroup $H\le\F_{k}$
is a directed, pointed, edge-labeled graph. This graph is called {\em
the core-graph associated with $H$} and is denoted by $\G_{X}(H)$.
We recall the notion of $\overline{\G}_{X}(H)$ the Schreier (right)
coset graph of $H$ with respect to the basis $X$. This is a directed,
pointed and edge-labeled graph. Its vertex set is the set of all right
cosets of $H$ in $\F_{k}$, where the basepoint corresponds to the
trivial coset $H$. For every coset $Hw$ and every letter $x_{i}$
there is a directed $i$-edge (short for $x_{i}$-edge) going from
the vertex $Hw$ to the vertex $Hwx_{i}$.

The core graph $\G_{X}(H)$ is obtained from $\overline{\G}_{X}(H)$
by omitting all the vertices and edges of $\overline{\G}_{X}(H)$
which are never traced by a reduced (i.e., non-backtracking) path
that starts and ends at the basepoint. Stated informally, we omit
all (infinite) {}``hanging trees'' from $\overline{\G}_{X}(H)$.
To illustrate, Figure \ref{fig:coset_and_core_graphs} shows the graphs
$\overline{\G}_{X}(H)$ and $\G_{X}(H)$ for $H=\langle x_{1}x_{2}x_{1}^{-3},x_{1}^{\;2}x_{2}x_{1}^{-2}\rangle\leq\F_{2}$.

\begin{figure}[h]
\begin{centering}
\begin{center}
\xy
(30,60)*+{\otimes}="s0";
(30,40)*+{\bullet}="s1";%
(45,30)*+{\bullet}="s2";%
(30,20)*+{\bullet}="s3";%
{\ar^{1} "s0";"s1"};%
{\ar^{2} "s1";"s2"};%
{\ar^{1} "s3";"s2"};%
{\ar^{1} "s1";"s3"};%
{\ar@(dl,dr)_{2} "s3";"s3"};%
(10,40)*+{\bullet}="t0";%
(60,40)*+{\bullet}="t1";%
(60,20)*+{\bullet}="t2";%
{\ar^{2} "t0";"s1"};%
{\ar^{1} "s2";"t1"};%
{\ar^{2} "s2";"t2"};%
{\ar@{..} "s0"; (20,60)*{}};%
{\ar@{..} "s0"; (30,70)*{}};%
{\ar@{..} "s0"; (40,60)*{}};%
{\ar@{..} "t0"; (10,50)*{}};%
{\ar@{..} "t0"; (0,40)*{}};%
{\ar@{..} "t0"; (10,30)*{}};%
{\ar@{..} "t1"; (52.5,45)*{}};%
{\ar@{..} "t1"; (67.5,45)*{}};%
{\ar@{..} "t1"; (67.5,35)*{}};%
{\ar@{..} "t2"; (67.5,25)*{}};%
{\ar@{..} "t2"; (67.5,15)*{}};%
{\ar@{..} "t2"; (52.5,15)*{}};%
{\ar@{=>} (76,40)*{}; (90,40)*{}};
(100,60)*+{\otimes}="m0";%
(100,40)*+{\bullet}="m1";%
(115,30)*+{\bullet}="m2";%
(100,20)*+{\bullet}="m3";%
{\ar^{1} "m0";"m1"};%
{\ar^{2} "m1";"m2"};%
{\ar^{1} "m3";"m2"};%
{\ar^{1} "m1";"m3"};%
{\ar@(dl,dr)_{2} "m3";"m3"};%
\endxy
\par\end{center}
\par\end{centering}

\caption{$\overline{\G}_{X}\left(H\right)$ and $\G_{X}\left(H\right)$ for
$H=\langle x_{1}x_{2}x_{1}^{-3},x_{1}^{\;2}x_{2}x_{1}^{-2}\rangle\leq\F_{2}$.
The Schreier coset graph $\overline{\G}_{X}\left(H\right)$ is the
infinite graph on the left (the dotted lines represent infinite $4$-regular
trees). The basepoint {}``$\otimes$'' corresponds to the trivial
coset $H$, the vertex below it corresponds to the coset $Hx_{1}$,
the one further down corresponds to $Hx_{1}^{\;2}=Hx_{1}x_{2}x_{1}^{-1}$,
etc. The core graph $\G_{X}\left(H\right)$ is the finite graph on
the right, which is obtained from $\overline{\G}_{X}\left(H\right)$
by omitting all vertices and edges that are not traced by reduced
closed paths around the basepoint.}

\label{fig:coset_and_core_graphs} 
\end{figure}

Note that the graph $\overline{\G}_{X}(H)$ is $2k$-regular: Every
vertex has exactly one outgoing $j$-edge and one incoming $j$-edge
for every $1\le j\le k$. Every vertex of $\G_{X}(H)$ has \emph{at
most} one outgoing $j$-edge, and {\em at most} one incoming $j$-edge
for every $1\le j\le k$.

It is an easy observation that 
\[
\pi_{1}(\overline{\G}_{X}(H))=\pi_{1}(\G_{X}(H))\stackrel{canonically}{\cong}H
\]
 where the canonical isomorphism is given by associating words in
$\F_{k}$ to paths in the coset graph and in the core graph: We traverse
the path by following the labels of outgoing edges. For instance,
the path (from left to right) 
\[
\xy(0,0)*+{\bullet}="s0";(18,0)*+{\bullet}="s1";(36,0)*+{\bullet}="s2";(54,0)*+{\bullet}="s3";(72,0)*+{\bullet}="s4";(90,0)*+{\bullet}="s5";(108,0)*+{\bullet}="s6";(126,0)*+{\bullet}="s7";{\ar^{2}"s0";"s1"};{\ar^{2}"s1";"s2"};{\ar^{1}"s2";"s3"};{\ar_{2}"s4";"s3"};{\ar^{3}"s4";"s5"};{\ar^{2}"s5";"s6"};{\ar_{1}"s7";"s6"};\endxy
\]
 corresponds to the word $x_{2}^{\;2}x_{1}x_{2}^{-1}x_{3}x_{2}x_{1}^{-1}$.
(See also \cite{MVW07}, where this fact appears in a slightly different
language).

Core graphs were introduced by Stallings \cite{Sta83}. Our definition
is slightly different, in that we allow the basepoint to have degree
one.

In fact, a {}``tail'' in $\G_{X}(H)$, i.e., a path to the basepoint
can be eliminated by replacing $H$ by an appropriate conjugate. However,
we find it unnecessary and less elegant for our needs.

We now list some properties of core graph, most of which are proved
in at least one of \cite{Sta83,KM02,MVW07}. The remaining ones are
easy observations.
\begin{claim}
\label{clm:core-graphs-properties} Let $H$ be a subgroup of $\F_{k}$
with an associated core graph $\G=\G_{X}(H)$. The Euler Characteristic
of a graph, denoted $\chi(\cdot)$ is the number of vertices minus
the number of edges. Finally, $rk(H)$ denotes the rank of the group
$H$.\end{claim}
\begin{enumerate}
\item $rk(H)<\infty\Leftrightarrow\G$ is finite 
\item $rk(H)=1-\chi(\G)$ 
\item Let $\Lambda$ be a finite, pointed, directed graph with edges labeled
by $\{x_{1},\ldots,x_{k}\}$. Then $\Lambda$ is a core graph (corresponding
to some $J\le\F_{k}$) if and only if $\Lambda$ satisfies the following
three properties: 

\begin{enumerate}
\item $\Lambda$ is connected 
\item With the possible exception of the basepoint, every vertex has degree
at least $2$. 
\item For every $1\le j\le k$, no two $j$-edges share the same origin
nor the same terminus. 
\end{enumerate}
\item There is a one-to-one correspondence between subgroups of $\F_{k}$
and core graphs. 
\item There is a one-to-one correspondence between subgroups of $\F_{k}$
of finite rank and finite core graphs. 
\end{enumerate}
In Appendix \ref{sec:Stallings-Folding-Algorithm} we present a well
known algorithm, based on Stallings' foldings, to obtain the core
graph of every $H\fg\F_{K}$ given some finite generating set for
$H$.

\subsection{Morphisms of Core Graphs}

\label{sbs:morphisms-of-core-graphs} In our framework, a morphism
between two core-graphs $\G_{1}$ and $\G_{2}$ is a map that sends
vertices to vertices and edges to edges, and preserves the structure
of the graphs. Namely, it preserves the incidence relations, sends
the basepoint to the basepoint, and preserves the directions and labels
of the edges.

As in Claim \ref{clm:core-graphs-properties}, the proofs of the following
properties are either easy variations on proofs in \cite{Sta83,KM02,MVW07}
or just easy observations:
\begin{claim}
\label{clm:morphism-properties} Let $H_{1},H_{2}\le\F_{k}$ be subgroups,
and $\G_{1},\G_{2}$ be the corresponding core graphs. Then \end{claim}
\begin{enumerate}
\item A morphism $\eta:\G_{1}\to\G_{2}$ exists $\Leftrightarrow H_{1}\leq H_{2}$,
\\
 and in this case, $\eta_{*}:\pi_{1}(\G_{1})\to\pi_{1}(\G_{2})$ is
injective. 
\item If a morphism exists, it is unique. 
\item Every morphism in an immersion (locally injective at the vertices). 
\end{enumerate}

\subsection{Quotients of Core Graphs}

\label{sbs:quotients-of-core-graphs} With core-graph morphisms at
hand, we can define the following rather natural relation between
core-graphs.
\begin{defn}
\label{def:quotient} Let $\G_{1},\G_{2}$ be core graphs and $H_{1},H_{2}\le\F_{k}$
the corresponding subgroups. We say that $\G_{1}$ \textbf{covers}
$\G_{2}$ or that $\G_{2}$ is a \textbf{quotient} of $\G_{1}$ if
there is a surjective morphism $\eta:\G_{1}\twoheadrightarrow\G_{2}$.
We also say in this case that $H_{1}$ \textbf{covers} $H_{2}$, and
denote $\G_{1}\twoheadrightarrow\G_{2}$ or $H_{1}\covers H_{2}$.
\end{defn}
By {}``surjective'' we mean surjective on both the vertices and
the edges. Note that we use the term {}``covers'' even though this
is \emph{not} a covering map in general (the morphism from $\G_{1}$
to $\G_{2}$ is always locally injective at the vertices, but not
necessarily locally bijective).

\begin{figure}[h]
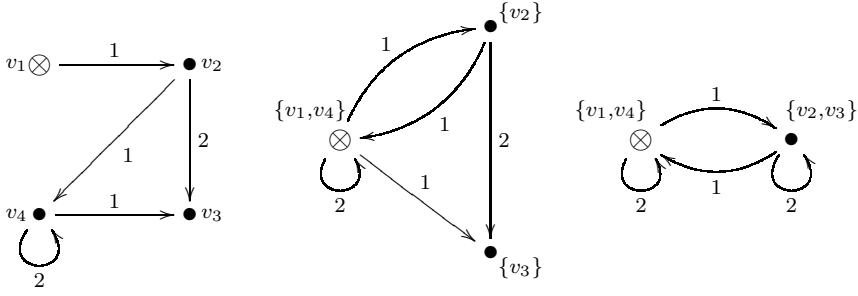

\noindent \begin{centering}
\begin{minipage}[t]{0.9\columnwidth}%
\noindent \begin{center}
\begin{center}
\xy 
(0,35)*+{\otimes}="m0"+(-3,0)*{\scriptstyle v_1};%
(20,35)*+{\bullet}="m1"+(3,0)*{\scriptstyle v_2};%
(20,15)*+{\bullet}="m2"+(3,0)*{\scriptstyle v_3};%
(0,15)*+{\bullet}="m3"+(-3,0)*{\scriptstyle v_4};%
{\ar^{1} "m0";"m1"};%
{\ar^{2} "m1";"m2"};%
{\ar^{1} "m3";"m2"};%
{\ar^{1} "m1";"m3"};%
{\ar@(dl,dr)_{2} "m3";"m3"};%
(40,25)*+{\otimes}="t0"+(-4,4)*{\scriptstyle \{v_1,v_4\}};%
(60,40)*+{\bullet}="t1"+(4,2)*{\scriptstyle \{v_2\}};%
(60,10)*+{\bullet}="t2"+(4,-2)*{\scriptstyle \{v_3\}};%
{\ar@/^1pc/^{1} "t0";"t1"};%
{\ar^{2} "t1";"t2"};%
{\ar^{1} "t0";"t2"};%
{\ar@/^1pc/^{1} "t1";"t0"};%
{\ar@(dl,dr)_{2} "t0";"t0"};%
(80,25)*+{\otimes}="t0"+(-4,4)*{\scriptstyle \{v_1,v_4\}};%
(100,25)*+{\bullet}="t1"+(4,4)*{\scriptstyle \{v_2,v_3\}};%
{\ar@/^1pc/^{1} "t0";"t1"};%
{\ar@/^1pc/^{1} "t1";"t0"};%
{\ar@(dl,dr)_{2} "t0";"t0"};%
{\ar@(dl,dr)_{2} "t1";"t1"};%
\endxy 
\par\end{center}
\par\end{center}%
\end{minipage}
\par\end{centering}

\caption{\label{fig:quotient-graph}The left graph is the core graph $\G_{X}(H)$
of $H=\left\langle x_{1}x_{2}x_{1}^{-3},x_{1}^{\;2}x_{2}x_{1}^{-2}\right\rangle \leq\F_{2}$.
Its vertices are denoted $v_{1},\ldots,v_{4}$. The graph in the middle
is the quotient corresponding to the partition $P=\{\{v_{1},v_{4}\},\{v_{2}\},\{v_{3}\}\}$.
This is not a core graph as there are two different $1$-edges originating
at $\{v_{1},v_{4}\}$. In order to obtain a core quotient-graph, we
use the folding process illustrated in Appendix \ref{sec:Stallings-Folding-Algorithm}.
The resulting core graph is on the right, corresponding to the partition
$P'=\{\{v_{1},v_{4}\},\{v_{2},v_{3}\}\}$.}
\end{figure}

For instance, $H=\langle x_{1}x_{2}x_{1}^{-3},x_{1}^{\;2}x_{2}x_{1}^{-2}\rangle\le\F_{k}$
covers the group \\
$J=\langle x_{2},x_{1}^{\;2},x_{1}x_{2}x_{1}\rangle$, the corresponding
core graphs of which are the leftmost and rightmost graphs in Figure
\ref{fig:quotient-graph}. As another example, every core graph $\G$
that contains edges of all labels covers the wedge graph $\Delta_{k}$.

We already know (Claim \ref{clm:morphism-properties}) that if $H_{1}\covers H_{2}$
then, in particular, $H_{1}\le H_{2}$. However, the converse is incorrect.
For example, the group\\
 \mbox{$K=\langle x_{1}x_{2}x_{1}^{-3},x_{1}^{\;2}x_{2}x_{1}^{-2},x_{2}\rangle$}
contains \mbox{$H=\langle x_{1}x_{2}x_{1}^{-3},x_{1}^{\;2}x_{2}x_{1}^{-2}\rangle$}
(we simply added $x_{2}$ as a third generator), yet $K$ is not a
quotient of $H$: the morphism $\eta:\G_{X}(H)\to\G_{X}(K)$ does
not contain the $2$-loop at the basepoint of $\G_{X}(K)$ in its
image.

Note also that the relation $H_{1}\covers H_{2}$ \emph{depends on
the given generating set $X$ of $\F_{k}$}. For example, if $H=\langle x_{1}x_{2}\rangle$
then $H\covers\langle x_{1},x_{2}\rangle=F_{2}$. However, $x_{1}x_{2}$
is primitive and could be taken as part of the original basis of $\F_{2}$.
In that case, the core graph of $H$ would consist of a single vertex
and single loop and would have no quotients except for itself.

It is also interesting to note that every quotient of the core-graph
$\G$ corresponds to some partition of $V(\G)$ (the partition determined
by the fibers of the morphism). We can simply draw a new graph with
a vertex for each block in the partition, and a $j$-edge from block
$b_{1}$ to block $b_{2}$ whenever there is some $j$-edge $(v_{1},v_{2})$
in $\G_{1}$ with $v_{1}\in b_{1},v_{2}\in b_{2}$. However, not every
partition of $V(\G)$ corresponds to a quotient core-graph: In the
resulting graph two distinct $j$-edges may have the same origin or
the same terminus. Note that even if a partition $P$ of $V(\G)$
yields a quotient which is not a core-graph, this can be remedied.
We can activate the folding process exemplified in Appendix \ref{sec:Stallings-Folding-Algorithm}
and obtain a core graph. The resulting partition $P'$ of $V(\G)$
is the finest partition which yields a quotient core-graph and which
is still coarser than $P$. We illustrate this in Figure \ref{fig:quotient-graph}.

\begin{lem}
\label{lem:finite-number-of-quotients} Every finite core-graph has
a finite number of quotients.\\
 Equivalently, every $H\fg\F_{k}$ covers a finite number of other
subgroups. \end{lem}
\begin{proof}
 The number of quotients of $\G$ is bounded from above by the number
of partitions of $V(\G)$. 
\end{proof}
Following the notations in \cite{MVW07}, we call the set of $X$-quotients
of $H$ the \emph{$X$-fringe} of $H$ and denote $\O_{X}(H)$. Namely,

\begin{equation}
\O_{X}(H):=\{\G_{X}(J)~|~H\covers J\}\label{eq:def_of_O_H}
\end{equation}
 Lemma \ref{lem:finite-number-of-quotients} states in this terminology
that for every $H\fg\F_{k}$ (and every basis $X$), $|\O_{X}(H)|<\infty$.

Before describing our new perspective on core graphs, we remind some
useful facts about free factors in free groups:
\begin{claim}
\label{clm:free_factors} Let $H,J,K\le\F_{k}$. Then, \end{claim}
\begin{enumerate}
\item Free factorness is transitive: If $H\ff J\ff K$ then $H\ff K$. 
\item If $\eta:\G_{X}(H)\hookrightarrow\G_{X}(J)$ is an embedding then
$H\ff J$. 
\item If $H\ff J$ then $H$ is a free factor in any subgroup $H\le M\le J$
in between. \end{enumerate}
\begin{proof}
The first and second claims are immediate. We give a {}``graph-theoretic''
proof for the third one. Assume that $H\ff J$, and let $Y$ be a
basis of $J$ extending some basis of $H$. In particular, $\G_{Y}(H)$
and $\G_{Y}(J)$ are both bouquets, consisting of a single vertex
and $rk(H)$ (resp. $rk(J)$) loops. Now, for every $H\le M\le J$,
consider the morphism $\eta:\G_{Y}(H)\to\G_{Y}(M)$. It is easy to
see that a core-graph-morphism of a bouquet must be an embedding.
Thus, by the second claim, $H\ff M$.
\end{proof}

\section{Immediate Quotients and the DAG of Core Graphs}

\label{sec:dag} The quotient relation yields a partial order on the
set of core graphs. But we are interested in a special case which
we call \emph{immediate quotients}. This relation allows us to build
the aforementioned DAG (directed acyclic graph) of all (core graphs
corresponding to) finite rank subgroups of $\F_{k}$.

Let $\G$ be a core graph, and let $P$ be a partition of $V(\G)$.
Let $\Delta$ be the quotient core graph we obtain from $P$ by the
folding process described in Figures \ref{fig:folding_process} and
\ref{fig:quotient-graph}. We say that $\Delta$ is \emph{generated}
from $\G$ by $P$. We are interested in the case where $P$ identifies
only a single pair of vertices:
\begin{defn}
\label{def:imme-quot} Let $\G$ be a core graph and let $P$ be a
partition of $V(\G)$ in which all parts consist of a single vertex
with a single exceptional part that contains two vertices. Let $\Delta$
be the core graph generated by $P$. We then say that $\Delta$ is
an \textbf{immediate quotient} of $\G$. 
\end{defn}
Alternatively we say that $\Delta$ is {\em generated by merging
a single pair} of vertices of $\G$. For instance, the rightmost
core graph in Figure \ref{fig:quotient-graph} is an immediate quotient
of the leftmost core graph.

The relation of immediate quotients has an interesting interpretation
for the associated free groups. Let $H,J\le\F_{k}$ be free groups
and $\G=\G_{X}(H),\Delta=\G_{X}(J)$ their core graph, and assume
$\Delta$ is an immediate quotient of $\G$ obtained by identifying
the vertices $u,v\in V(\G)$. Now let $p_{u},p_{v}\in\F_{k}$ be words
corresponding to some paths in $\G$ from the basepoint to $u$ and
$v$ respectively. It is not hard to see that identifying $u$ and
$v$ has the same effect as adding the word $w=p_{u}p_{v}^{-1}$ to
$H$ and considering the generated group. Namely, $J=\langle H,w\rangle$.

\begin{center}
\begin{center}
\xy 
(0,0)*{\otimes}="bp"+(-8,8)*{\G};%
(5,0)*\xycircle(30,15){--};
(24,4)*{\bullet}="u"  +(2,2)*{u};%
(20,-6)*{\bullet}="v" +(2,-2)*{v};%
"bp"; "u" **\crv{(10,20) & (15,-20)} ?(.43)*\dir{>>}; 
(6,8)*{\scriptstyle p_u}; "bp"; "v" **\crv{(5,-10) & (15,10)} ?(.43)*\dir{>>}; 
(6,-6)*{\scriptstyle p_v}; (-60,0)*+{}="dummy";
(0,17)*+{}="dummy";
(0,-17)*+{}="dummy";
\endxy
\par\end{center}
\par\end{center}

Based on the relation of immediate quotients we consider the DAG $\D_{k}$.
The set of vertices of this graph consists of all finite core graphs
with edges labeled by $1,\ldots,k$, and its directed edges connect
every core graph to its immediate quotients. Every fixed ordered basis
of $\F_{k}$ $X=\{x_{1},\ldots,x_{k}\}$, determines a one-to-one
correspondence between the vertices of this graph and all finite rank
subgroups of $\F_{k}$.

\begin{figure}[h]
\begin{centering}
\includegraphics{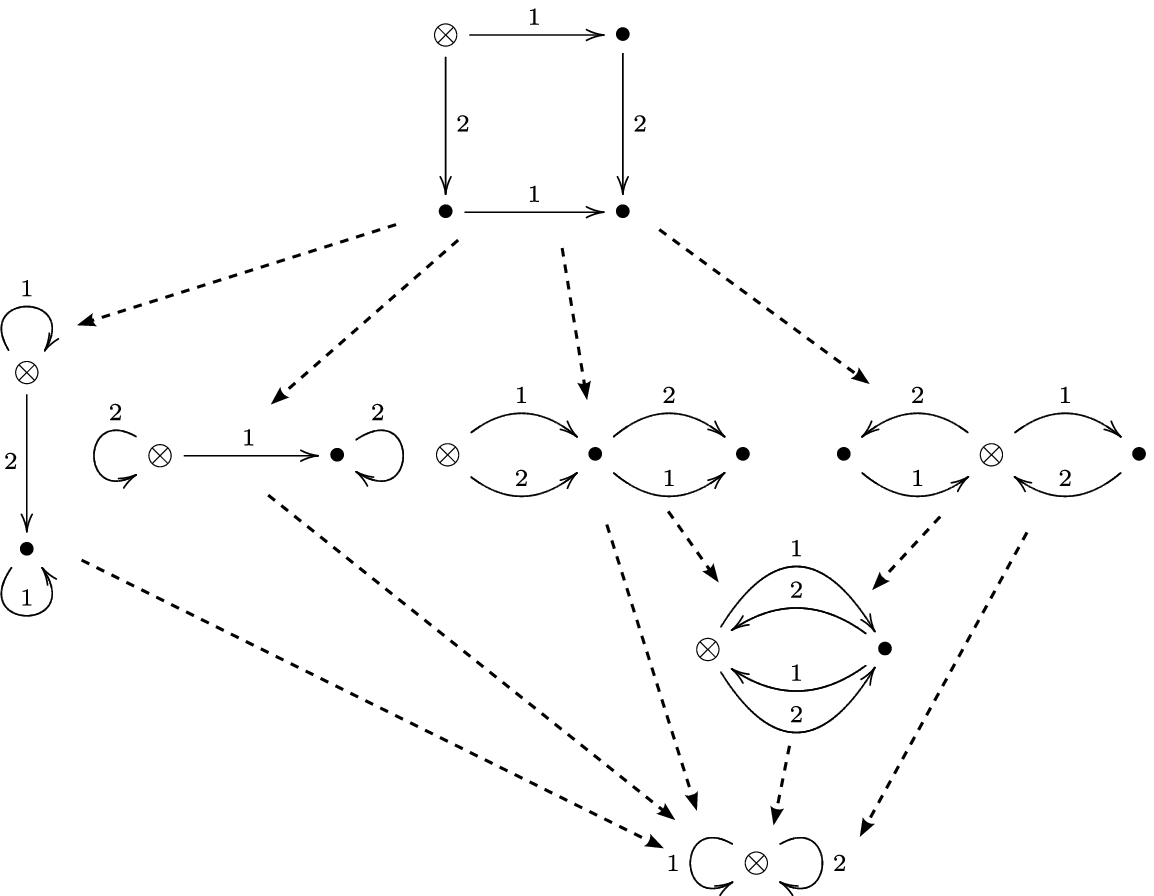}
\par\end{centering}

\caption{The subgraph of $\D_{k}$ induced by $\O_{X}(H)$, that is, all quotients
of the core graph $\G=\G_{X}\left(H\right)$, for $H=\langle x_{1}x_{2}x_{1}^{-1}x_{2}^{-1}\rangle$.
The dashed arrows denote immediate quotients, i.e.\ quotients generated
by merging a single pair of vertices. $\G$ has exactly seven quotients:
itself, four immediate quotients, and two quotients at distance $2$.}

\label{fig:lattice} 
\end{figure}

As before, we fix an ordered basis $X$. For any $H\fg\F_{k}$, the
subgraph of $\D_{k}$ of the descendants of $\G_{X}(H)$ consists
of all quotients of $\G_{X}(H)$, that is of all elements of the $X$-fringe
$\O_{X}(H)$. By Lemma \ref{lem:finite-number-of-quotients}, this
subgraph is finite. In Figure \ref{fig:lattice} we draw the subgraph
of $\D_{k}$ consisting of all quotients of $\G_{X}(H)$ when $H=\langle x_{1}x_{2}x_{1}^{-1}x_{2}^{-1}\rangle$.
The edges of this subgraph (i.e. immediate quotients) are denoted
by the broken arrows in the figure.

This yields the aforementioned distance function between a finite
core graph and a quotient of it:
\begin{defn}
\label{def:distance} Let $H_{1},H_{2}\fg\F_{k}$ be finite rank subgroups
such that $H_{1}\covers H_{2}$, and let $\G_{1}=\G_{X}(H_{1}),\G_{2}=\G_{X}(H_{2})$
be the corresponding core graphs. We define the distance between $H_{1}$
and $H_{2}$, denoted $\rho_{X}(H_{1},H_{2})$ or $\rho(\G_{1},\G_{2})$,
to be the shortest length of a directed path from $\G_{1}$ to $\G_{2}$
in $\D_{k}$. 
\end{defn}
In other words, $\rho_{X}(H_{1},H_{2})$ is the length of the shortest
series of immediate quotients that yields $\G_{2}$ from $\G_{1}$.
Equivalently, it is the minimal number of pairs of vertices that need
to be identified in $\G_{1}$ in order to obtain $\G_{2}$ (via the
folding process). For example, if $\G_{2}$ is an immediate quotient
of $\G_{1}$ then $\rho_{X}(H_{1},H_{2})=\rho(\G_{1},\G_{2})=1$.
For $H=\langle x_{1}x_{2}x_{1}^{-1}x_{2}^{-1}\rangle$, $\G_{X}(H)$
has four quotients at distance $1$ and two at distance $2$ (see
Figure \ref{fig:lattice}).

As aforementioned, by merging a single pair of vertices of $\G_{X}(H)$
(and then folding) we obtain the core graph of a subgroup $J$ obtained
from $H$ by adding some single generator. Thus, by taking an immediate
quotient, the rank of the associated subgroup increases at most by
1 (in fact, it may also stay unchanged or even decrease). This implies
that whenever $H\covers J$:

\begin{equation}
rk(J)-rk(H)~~\le~~\rho(H,J)\label{eq:rho>=00003Drk-rk}
\end{equation}

It is not hard to bound the distance from above as well:
\begin{lem}
\label{lem:bounds_for_rho} Let $H,J\fg\F_{k}$ such that $H\covers J$.
Then 
\[
rk(J)-rk(H)~~\le~~\rho_{X}(H,J)~~\le~~rk(J)
\]

\end{lem}
We postpone the proof of the upper bound to Appendix \ref{sec:The-Proof-of-bounds-lemma}.
(In fact, this upper bound in not needed for the main results of this
paper. We give it anyway in order to have the full picture in mind.)
Theorem \ref{thr:rho=00003Drk-rk_iff_ff} then states that in the
same setting, the lower bound is attained iff $H$ is a free factor
of $J$. In fact one can visualize these results in the following
way. Consider a two dimensional table which contains all the elements
of the fringe $\O_{X}(H)$ (each quotient of $\G_{X}(H)$ contained
in some, not necessarily distinct, cell). The rows correspond to the
rank and are indexed $0,1,2,3,...$. The columns correspond to the
distance from $H$ and are also indexed $0,1,2,3,\ldots$. We then
put every quotient of $H$ in the suitable cell in the table. Let
$r=rk(H)$ denote the rank of $H$. Lemma \ref{lem:bounds_for_rho}
then says that the (finitely many) elements of $\O_{X}(H)$ are spread
across $r+1$ diagonals in the table: the main one and the $r$ diagonals
below it. Theorem \ref{thr:rho=00003Drk-rk_iff_ff} implies that within
$\O_{X}(H)$, $H$ is a free factor of exactly those $J$-s found
in the lowest of these $r+1$ diagonals. (In fact, Lemma \ref{lem:pi=00003Dpi'}
shows that $\pi(H)$ can also be read from this table: it equals the
rank of the upmost occupied cell in this table outside the free-factor-diagonal.)

\subsection{Proof of Theorem \ref{thr:rho=00003Drk-rk_iff_ff}}

\label{sbs:proof_rho=00003Drk-rk_iff_ff}

The main result of this section states that if $H\fg J\fg\F_{k}$
and $H\covers J$, then 
\begin{equation}
\rho_{X}(H,J)=rk(J)-rk(H)~~~\Longleftrightarrow~~~H\ff J.\label{eq:rho=00003Drk-rk_iff_ff}
\end{equation}

In fact, one of the implications is trivial. As mentioned above, merging
two vertices in $\G_{X}(H)$ is equivalent to adding some generator
to $H$. If we manage to obtain $\G_{X}(J)$ from $\G_{X}(H)$ by
$\mathrm{rk}(J)-\mathrm{rk}(H)$ merging steps, this means we can
obtain $J$ from $H$ by adding $rk(J)-rk(H)$ extra generators to
$H$, hence clearly $H\ff J$ (recall that by hopfianity of the free
group, every generating set of size $\mathrm{rk}\left(J\right)$ is
a basis of $J$, see e.g.~\cite[Chapter 2.29]{Bog08}). Thus, 
\begin{equation}
\rho_{X}(H,J)=rk(J)-rk(H)~~~\Longrightarrow~~~H\ff J\label{eq:easy-implication}
\end{equation}

The other implication is not trivial. Assume that $H\ff J$. Our goal
is to obtain $rk(J)-rk(H)$ complementary generators of $J$ from
$H$, so that each of them can be realized by merging a pair of vertices
in $\G_{X}(H)$. 

To this goal we introduce the notion of a {}``\emph{handle number}''
associated with a subgroup $M$ and a word $w\in\F_{k}$. (It also
depends on the fixed basis $X$ of $\F_{k}$). This number is defined
as follows. Let $\G=\G_{X}(M)$. Denote by $p_{w}$ the longest prefix
of $w$ that corresponds to some path from the basepoint of $\G$
(we trace the letters of $w$ along $\G$ until we get stuck). Likewise,
denote by $s_{w}$ the longest suffix of $w$ that ends at the basepoint
(here we trace $w^{-1}$ from the basepoint until we get stuck). If
$|p_{w}|+|s_{w}|<|w|$, then $w=p_{w}m_{w}s_{w}$ as a reduced expression
for some $1\neq m_{w}\in\F_{k}$. The handle number of $\left(M,w\right)$
is then
\[
h_{X}\left(M,w\right)=h\left(\Gamma,w\right)=\begin{cases}
\left|m_{w}\right| & |p_{w}|+|s_{w}|<|w|\\
0 & \mathrm{otherwise}
\end{cases}.
\]

\begin{claim}
\label{claim:1-handle}Assume that $w\notin M$ and let $N=\left\langle M,w\right\rangle $.
Then
\begin{enumerate}
\item $h_{X}\left(M,w\right)>0$ if and only if $\Gamma_{X}\left(M\right)$
is a (proper) subgraph of $\Gamma_{X}\left(N\right)$, and
\item $h_{X}\left(M,w\right)=0$ if and only if $\Gamma_{X}\left(N\right)$
is an immediate quotient of $\Gamma_{X}\left(M\right)$.
\end{enumerate}
\end{claim}
\begin{proof}
Assume first that $h_{X}\left(M,w\right)>0$. In the notations of
the previous paragraph, let $v(p_{w}),v(s_{w})$ be the end point
of the path corresponding to $p_{w}$ and the starting point of the
path corresponding to $s_{w}$. We can then add a {}``handle'' to
$\G$ in the form of a path corresponding to $m_{w}$ which starts
at $v(p_{w})$ and ends at $v(s_{w})$. (If $v(p_{w})=v(s_{w})$ this
handle looks like a balloon, possibly with a string.)

\begin{center}
\begin{center}
\xy (0,0)*{\otimes}="bp"+(-8,8)*{\G};%
(5,0)*\xycircle(20,20){--};
(22.32,10)*{\bullet}="s1"  +(4,4)*{v(p_w)};%
(22.32,-10)*{\bullet}="s2" +(4,-4)*{v(s_w)};%
"bp"; "s1" **\crv{(10,20) & (15,-20)} ?(.43)*\dir{>>}; 
(6,8)*{\scriptstyle p_w}; 
"s2"; "bp" **\crv{(20,20) & (5,-10)} ?(.33)*\dir{>>}; 
(6,-4)*{\scriptstyle s_w}; 
"s1";"s2" **\crv{(40,10) & (40,-10)} ?(.5)*\dir{>>}; 
(40,0)*{\scriptstyle m_w}; 
(-40,20)*+{}="dummy";
(0,22)*+{}="dummy";
(0,-22)*+{}="dummy";
\endxy 
\par\end{center}
\par\end{center}

The resulting graph is a core-graph (the edge conditions at $v(p_{w})$
and $v(s_{w})$ hold, by the maximality of $p_{w}$ and $s_{w}$),
and it corresponds to $N$. So we actually obtained $\Gamma_{X}\left(N\right)$.
It follows that $\G_{X}\left(M\right)$ is a proper subgraph of $\Gamma_{X}\left(N\right)$.
On the other hand, if $h_{X}\left(M,w\right)=0$, i.e.~if $|p_{w}|+|s_{w}|\geq|w|$,
we can find a pair of vertices in $\G$ whose merging adds $w$ to
$H$ as a complementary generator for $J$. (We may take $v(p_{w})$
together with the vertex on the path of $s_{w}$ at distance $|p_{w}|+|s_{w}|-|w|$
from $v(s_{w})$.) 

\end{proof}
The last claim shows, in particular, that if $N$ is obtained from
$M$ by adding a single complementary generator, then either $\Gamma_{X}\left(N\right)$
is an immediate quotient or it contains $\Gamma_{X}\left(M\right)$
as a proper subgraph. This already proves Theorem \ref{thr:rho=00003Drk-rk_iff_ff}
for the case $\mathrm{rk}\left(J\right)-\mathrm{rk}\left(H\right)=1$:
if $H\covers J$, we are clearly in the second case of Claim \ref{claim:1-handle},
i.e.~$J$ is an immediate quotient of $H$.

We proceed by defining handle numbers for a subgroup $M\fg\F_{k}$
and an ordered set of words $w_{1},\ldots,w_{t}\in\F_{k}$. Let $N=\left\langle M,w_{1},\ldots,w_{t}\right\rangle $
and $\Delta=\Gamma_{X}\left(N\right)$. Let in addition $N_{i}=\langle M,w_{1},\ldots,w_{i}\rangle$
and $\G_{i}=\G_{X}(N_{i})$. We obtain a series of subgroups 
\[
M=N_{0}\leq N_{1}\leq\ldots\leq N_{t}=N,
\]
 and a series of graphs $\G=\G_{0},\G_{1},\ldots,\G_{t}=\Delta$.
We denote by $h_{X}\left(M,w_{1},\ldots,w_{t}\right)$ the $t$-tuple
of the following handle numbers:

\[
h_{X}\left(M,w_{1},\ldots,w_{t}\right)\stackrel{\mathrm{def}}{=}\left(h\left(\Gamma_{0},w_{1}\right),\, h\left(\Gamma_{1},w_{2}\right),\,\ldots,h_{X}\left(\Gamma_{t-1},w_{t}\right)\right).
\]

Let us focus now on the case where t is the cardinality of the smallest
set $S\subseteq\F_{k}$ such that $N=\left\langle M,S\right\rangle $.
The following lemma characterizes $t$-tuples of words for which the
$t$-tuple of handle-numbers is lexicographically minimal. It is the
crux of the proof of Theorem \ref{thr:rho=00003Drk-rk_iff_ff}.
\begin{lem}
\label{prop:handle_numbers}In the above notations, let $ $$\left(w_{1},\ldots,w_{t}\right)$
be an ordered set of complementary generators such that the tuple
of handle numbers $h_{X}\left(M,w_{1},\ldots,w_{t}\right)$ is lexicographically
minimal. Then the zeros in $h_{X}\left(M,w_{1},\ldots,w_{t}\right)$
form a prefix of the tuple.
\end{lem}
Namely, there is no zero handle-number that follows a positive handle-number.
\begin{proof}
It is enough to prove the claim for pairs of words (i.e.~for $t=2$),
the general case following immediately. Assume then that $N=\left\langle M,w_{1},w_{2}\right\rangle $,
that $2$ is the minimal number of complementary generators of $N$
given $M$, and that $h_{X}\left(M,w_{1},w_{2}\right)$ is lexicographically
minimal. In the above notation, assume to the contrary that $h(\G_{0},w_{1})>0$
and $h(\G_{1},w_{2})=0$. Let $m_{1}=m_{w_{1}}$ denote the handle
of $w_{1}$ in $\Gamma_{0}$. Thus $\G_{1}$ was obtained from $\G_{0}$
by adding a handle (or a balloon) representing $m_{1}$. The word
$w_{2}$ can be expressed as $w_{2}=ps$ so that there is a path corresponding
to $p$ in $\Gamma_{1}$, emanating from the basepoint and ending
at some vertex $u$, and there is a path $s$ to the basepoint from
a vertex $v$. (Clearly, $u\ne v$ for otherwise $w_{2}\in N_{1}$
contradicting the minimality of $t=2$.) Thus $\G_{2}$ is attained
from $\G_{1}$ by identifying the vertices $u$ and $v$. By possibly
multiplying $w_{2}$ from the left by a suitable element of $N_{1}$,
we can assume that $p$ does not traverse the handle $m_{1}$ {}``more
than necessary''. Namely, if $u$ does not lie on $m_{1}$, then
$p$ avoids $m_{1}$, and if it does lie on $m_{1}$, then only the
final segment of $p$ traverses $m_{1}$ till $u$. The same holds
for $s$ and $v$ (with right multiplication).

The argument splits into three possible cases. 
\begin{itemize}
\item If both $u,v$ belong already to $\G_{0}$, then $h(\G_{0},w_{2})=0$.
In this case we can switch between $w_{2}$ and $w_{1}$ to lexicographically
reduce the sequence of handle numbers, contrary to our assumption. 
\item Consider next the case where, say, $v\in V(\G_{0})$ but $u\in V(\G_{1})\setminus V(\G_{0})$,
i.e, $u$ resides on the handle $m_{1}$. Then, the handle needed
in order to add $w_{2}$ to $\G_{0}$ is strictly shorter than $h(\G_{0},w_{1})=\left|m_{1}\right|$.
Again, by switching $w_{2}$ with $w_{1}$ the sequence of handle
numbers goes down lexicographically - a contradiction.
\end{itemize}

\begin{center}
\begin{center}
\xy 
(0,0)*{\otimes}="bp"+(-6,8)*{\G_{0}};%
(5,0)*\xycircle(20,20){--};
(22.32,10)*{}="s1";%
(22.32,-10)*{}="s2";%
"s1";"s2" **\crv{(40,10) & (40,-10)}  ?(.2)*\dir{>>}; \POS?(.40)*{\bullet}="u"+(3,3)*{u}; 
(40,0)*{\scriptstyle m_{1}}; 
(15,-17.32)*{\bullet}="v"+(3,-3)*{v}; 
"bp"; "s1" **\crv{(10,20) & (15,-20)} ?(.43)*\dir{>>}; 
(30,12)*{\scriptstyle p}; 
"v"; "bp" **\crv{(20,20) & (5,-10)} ?(.33)*\dir{>>}; 
(6,-4)*{\scriptstyle s}; 
(-40,20)*+{}="dummy";
(0,22)*+{}="dummy";
(0,-22)*+{}="dummy";
\endxy
\par\end{center}
\par\end{center}
\begin{itemize}
\item In the final case that should be considered both $u$ and $v$ are
on the handle $m_{1}$. I.e. $u,v\in V(\G_{1})\setminus V(\G_{0})$.
Assume w.l.o.g.~that when tracing the path of $m_{1}$, $u$ precedes
$v$. As before we can premultiply and postmultiply $w_{2}$ by suitable
elements of $N_{1}$ to guarantee the following: The path $p$, from
the basepoint of $\G_{1}$ to $u$, goes through $\Gamma_{0}$ and
then traverses a prefix of $m_{1}$ until reaching $u$, and the path
$s$ from $v$ to the basepoint traces a suffix of $m_{1}$ and then
goes only through $\Gamma_{0}$. Again $h(\G_{0},w_{2})<h(\G_{0},w_{1})$,
so that switching $w_{2}$ with $w_{1}$ lexicographically reduces
the sequence of handle numbers. (A similar argument works in the case
$m_{1}$ constitutes a balloon.)
\end{itemize}

\begin{center}
\begin{center}
\xy 
(0,0)*{\otimes}="bp"+(-6,8)*{\G_{0}};%
(5,0)*\xycircle(20,20){--};
(22.32,10)*{}="s1";%
(22.32,-10)*{}="s2";%
"s1";"s2" **\crv{(40,10) & (40,-10)}  ?(.2)*\dir{>>}; ?(.9)*\dir{>>}; \POS?(.40)*{\bullet}="u"+(3,3)*{u}; 
\POS?(.75)*{\bullet}="v"+(2,-3)*{v}; 
(40,0)*{\scriptstyle m_{1}}; 
"bp"; "s1" **\crv{(10,20) & (15,-20)} ?(.43)*\dir{>>}; 
(30,12)*{\scriptstyle p}; 
"s2"; "bp" **\crv{(20,20) & (5,-10)} ?(.33)*\dir{>>}; 
(28,-12)*{\scriptstyle s}; 
(-40,20)*+{}="dummy";
(0,22)*+{}="dummy";
(0,-22)*+{}="dummy";
\endxy
\par\end{center}
\par\end{center}

\end{proof}
Theorem \ref{thr:rho=00003Drk-rk_iff_ff} follows easily from this
lemma. Indeed, assume that $H\ff J$ and that $H\covers J$. Let $t=\mathrm{rk}\left(J\right)-\mathrm{rk}\left(H\right)$
denote the difference in ranks, so that $t$ is the smallest number
of complementary generators needed to obtain $J$ given $H$. Let
$\left(w_{1},\ldots,w_{t}\right)$ be an ordered set of complementary
generators so that $h_{X}\left(H,w_{1},\ldots,w_{t}\right)$ is lexicographically
minimal. Similarly to the notations above, let $J_{i}=\left\langle H,w_{1},\ldots,w_{i}\right\rangle $
and $\Gamma_{i}=\Gamma_{X}\left(J_{i}\right)$. 

By the lemma, there is some $0\leq q\leq t$ so that $h\left(\Gamma_{0},w_{1}\right)=\ldots=h\left(\Gamma_{q-1},w_{q}\right)=0$
whereas $h(\Gamma_{q},w_{q+1}),\ldots,h\left(\Gamma_{t-1},w_{t}\right)$
are all positive. By Claim \ref{claim:1-handle} it follows that $\Gamma_{i}$
is an immediate quotient of $\Gamma_{i-1}$ for $1\leq i\leq q$,
and therefore $J_{q}\in\O_{X}\left(H\right)$ and $\rho_{X}\left(H,J_{q}\right)=q$.
(This in fact shows that $\rho_{X}(H,J_{q})\leq rk(J_{q})-rk(H)$,
and the equality follows from Lemma \ref{lem:bounds_for_rho}).

Using Claim \ref{claim:1-handle} again, we see that $\Gamma_{i}$
is a proper subgraph of $\Gamma_{i+1}$ for $q\leq i\leq t-1$. So
that $\Gamma_{q}$ is a subgraph of $\Gamma_{t}=\Gamma_{X}\left(J\right)$.
But then the image of the graph morphism $\eta:\Gamma_{X}\left(H\right)\to\Gamma_{X}\left(J\right)$
is clearly the subgraph $\Gamma_{q}$. If $q<t$ this is a proper
subgraph, which contradicts the assumption $H\covers J$. Hence $q=t$
and $\rho_{X}\left(H,J\right)=t$, as required. Together with \prettyref{eq:easy-implication}
this completes the proof of Theorem \prettyref{thr:rho=00003Drk-rk_iff_ff}.
$\qed$

In fact, the same argument yields a more general result:
\begin{cor}
\label{cor:min_num_of_complemetary}Let $H\leq J\leq\F_{k}$ be f.g.~groups,
and let $t$ be the minimal number of complementary generators needed
to obtain $J$ from $H$. Then $t$ is computable as follows. Let
$\eta:\Gamma_{X}\left(H\right)\to\Gamma_{X}\left(J\right)$ be the
unique core-graph morphism, and let $M$ be the intermediate subgroup
corresponding to the image $\eta\left(\Gamma_{X}\left(H\right)\right)$.
Then,
\[
t=\rho_{X}\left(H,M\right)+\mathrm{rk}\left(J\right)-\mathrm{rk}\left(M\right).
\]
\end{cor}
\begin{proof}
In the notation of the last part of the proof of Theorem \prettyref{thr:rho=00003Drk-rk_iff_ff},
we see that $M=J_{q}\ff J$, and so $\rho_{X}\left(H,M\right)+\mathrm{rk\left(J\right)}-\mathrm{rk}\left(M\right)=\rho_{X}\left(H,J_{q}\right)+\left(\mathrm{t-q}\right)=t$.\end{proof}
\begin{rem}
\label{remark:pairs-of-vertices}Note that in the crucial arguments
of the proof of Theorem \ref{thr:rho=00003Drk-rk_iff_ff} we did not
use the fact that the groups where of finite rank. Indeed, this result
can be carefully generalized to subgroups of $\F_{k}$ of infinite
rank.
\end{rem}

\begin{rem}
Another way to interpret Theorem \prettyref{thr:rho=00003Drk-rk_iff_ff}
is by saying that if $H\ff J$ and $H\covers J$ with $t=\mathrm{rk\left(J\right)-rk\left(H\right)}$,
then there exists some set $\left\{ w'_{1},\ldots,w'_{t}\right\} $
of complementary generators such that each $w_{i}$ can be realized
by merging a pair of vertices in $\Gamma_{X}\left(H\right)$. To see
this, let $w_{1},\ldots,w_{t}$ be as in the proof above, so $w_{i}$
can be realized by merging a pair of vertices $u$ and $v$ in $\Gamma_{i-1}$.
Let $\eta_{i-1}:\Gamma_{X}\left(H\right)\to\Gamma_{i-1}$ be the surjective
morphism, and pick any vertices in the fibers $u'\in\eta^{-1}\left(u\right),\, v'\in\eta^{-1}\left(v\right)$.
Let $w_{i}'$ be some word corresponding to the merging of $u'$ and
$v'$ in $\Gamma_{X}\left(H\right)$. It is not hard to see that for
each $i$, $\left\langle H,w_{1},\ldots,w_{i}\right\rangle =\left\langle H,w_{1}',\ldots,w_{i}'\right\rangle $.
\end{rem}

\section{More on the Primitivity Rank}

\label{sec:prim-rank}

Recall Definitions \ref{def:primitivity-rank} and \ref{def:primitivity-rank-subgp}
where we defined $\pi(w)$, the primitivity rank of a word $w\in\F_{k}$,
and $\pi(H)$, the primitivity rank of $H\fg\F_{k}$. In this subsection
we prove some characteristics of this categorization of formal words,
and show it actually depends only on the quotients of the core graph
$\G_{X}(H)$ (or $\G_{X}(\la w\ra)$). The claims are stated for subgroups,
and can be easily interpreted for elements with the usual correspondence
between the element $w\ne1$ and the subgroup it generates $\la w\ra$.
We begin by characterizing the possible values of $\pi(H)$.
\begin{lem}
\label{lem:primitive-gets-infty} Let $H\fg\F_{k}$. Then 
\[
H\ff\F_{k}\Leftrightarrow\pi(H)=\infty.
\]
 \end{lem}
\begin{proof}
Recall that $\pi(H)$ is defined by the smallest rank of subgroups
of $\F_{k}$ where $H$ is contained but not as a free factor. If
$H$ is not a free factor of $\F_{k}$, then $\F_{k}$ itself is one
such subgroup so that $\pi(H)\leq k<\infty$. If $H\ff\F_{k}$, Claim
\ref{clm:free_factors} shows it is a free factor in every other subgroup
containing it. Thus, in this case $\pi(H)=\infty$. \end{proof}
\begin{cor}
\label{cor:image-of-pi} For every $H\fg\F_{k}$ 
\[
\pi(H)\in\{0,1,\ldots,k\}\cup\{\infty\}
\]

\end{cor}
In the definition of the primitivity rank of a subgroup $H$, we consider
all subgroups of $\F_{k}$ containing $H$ but not as a free factor.
It turns out it is enough to consider only subgroups of $\F_{k}$
that are covered by $H$, that is, groups whose associated core graphs
are in the $X$-fringe $\O_{X}(H)$.
\begin{lem}
\label{lem:pi=00003Dpi'} For every $H\fg\F_{k}$ 
\begin{equation}
\pi(H)=min\left\{ rk(J)~\Big|~\begin{gathered}H\covers J~~and\\
H\textrm{ is \textbf{not} a free factor of }J
\end{gathered}
\right\} \label{eq:pi=00003Dpi'}
\end{equation}
 Moreover, all $H$-critical subgroups of $\F_{k}$ are covered by
$H$. \end{lem}
\begin{proof}
Recall that $H$-critical subgroups of $\F_{k}$ are the subgroups
of smallest rank in which $H$ is not a free factor (so in particular
their rank is exactly $\pi(H)$). It is enough to show that every
$H$-critical subgroup has its associated core graph in the fringe
$\O_{X}(H)$.

Consider an $H$-critical subgroup $J\le\F_{k}$. This $J$ contains
$H$ but not as a free factor. By Claim \ref{clm:morphism-properties}
there exists a morphism $\eta:\G_{X}(H)\to\G_{X}(J)$. If $\eta$
is surjective then $H\covers J$ and $\G_{X}(J)\in\O_{X}(H)$. Otherwise,
consider $J'$, the group corresponding to the core graph $\eta(\G_{X}(H))$.
This graph, $\G_{X}(J')$, is a strict subgraph of $\G_{X}(J)$, and
so $J'\stackrel{*}{\lneqq}J$ (see Claim \ref{clm:free_factors}).
In particular $H\covers J'$ and $rk(J')<rk(J)$. It is impossible
that $H\ff J'$, because by transitivity this would yield that $H\ff J$.
Thus, $J'$ is a subgroup in which $H$ is a not free factor, and
of smaller rank than $J$. This contradicts the fact that $J$ is
$H$-critical. 
\end{proof}
We note that in the terminology of \cite{KM02,MVW07}, $H$-critical
subgroups are merely a special kind of {}``algebraic extensions''
of $H$. (An algebraic extension of $H$ is a group $J$ such that
for every $M$ with $H\le M\lvertneqq J$, $M$ is not a free factor
of $J$.) Specifically, $H$-critical subgroups are algebraic extensions
of $H$ of minimal rank, excluding $H$ itself. Our proof actually
shows the more general fact that all algebraic extensions of $H$
can be found in the fringe (this fact appears in \cite{KM02,MVW07}).

At this stage we can describe exactly how the primitivity rank of
a subgroup $H\fg\F_{k}$ can be computed. In fact, all algebraic extensions
and critical subgroups of $H$ can be immediately identified:
\begin{cor}
\label{cor:identify-alg-extensions}Consider the induced subgraph
of $\D_{k}$ consisting of all core graphs in $\O_{X}(H)$. Then,\end{cor}
\begin{itemize}
\item The algebraic extensions of $H$ are precisely the core graphs which
are not an immediate quotient of any other core graph of smaller rank. 
\item The $H$-critical subgroups are the algebraic extensions of smallest
rank, excluding $H$ itself, and $\pi\left(H\right)$ is their rank.\end{itemize}
\begin{proof}
The second statement follows from the discussion above and from definition
\ref{def:primitivity-rank-subgp}. The first statement holds trivially
for $H$ itself. If $J$ is a proper algebraic extension of $H$,
then by the proof of Lemma \ref{lem:pi=00003Dpi'}, $J\in\O_{X}\left(H\right)$.
If $\Gamma_{X}\left(J\right)$ is an immediate quotient of some $\Gamma_{X}\left(M\right)$
of smaller rank, where $M\in\O_{X}\left(H\right)$, then $H\leq M\lneq J$
and by (the easier implication of) Theorem \ref{thr:rho=00003Drk-rk_iff_ff}
we conclude $M\ff J$, a contradiction. 

On the other hand, if $J\in\O_{X}\left(H\right)$ is not an algebraic
extension of $H$, then there is some intermediate subgroup $L$ such
that $H\leq L\stackrel{*}{\lneqq}J$. We can assume $L\in\O_{X}\left(H\right)$
for otherwise it can be replaced with $L'$ corresponding to the image
of the morphism $\eta:\Gamma_{X}\left(H\right)\to\Gamma_{X}\left(L\right)$
(whence $L'\in\O_{X}\left(H\right)$ and $H\leq L'\ff L\stackrel{*}{\lneqq}J$).
From (the harder implication of) Theorem \prettyref{thr:rho=00003Drk-rk_iff_ff}
it follows that $\rho_{X}\left(L,J\right)=\mathrm{rk}\left(J\right)-\mathrm{rk}\left(L\right)$.
The prior-to-last element in a shortest path in $\D_{k}$ from $\Gamma_{X}\left(L\right)$
to $\Gamma_{X}\left(J\right)$ is then a proper free factor of $J$
at distance 1 that belongs to $\O_{X}\left(H\right)$.
\end{proof}
As an example, consider $H=\la x_{1}x_{2}x_{1}^{-1}x_{2}^{-1}\ra$.
The full lattice of groups in $\O_{X}(H)$ is given in Figure \ref{fig:lattice}.
There is one group of rank $1$ ($H$ itself), 5 of rank $2$ and
one of rank $3$. The only group in the lattice where $H$ in not
a free factor is the group $\la x_{1}x_{2}\ra=\F_{2}$, of rank $2$,
so $\pi(H)=2$. (And the set of algebraic extensions of $H$ is precisely
$\{H,\F_{2}\}$.)

\section{The Calculation of \texorpdfstring{$\phi$}{phi}}

\label{sec:phi} The proof of Proposition \ref{prop:phi=00003Dpi},
as well as the reasoning that underlies Conjecture \ref{conj:phi=00003Dpi},
are based on the fact that both $\phi(H)$ and $\pi(H)$ can be calculated
by analyzing $\O_{X}(H)$, the set of quotients of $\G_{X}(H)$. In
the previous section it was shown how $\pi(H)$ is determined from
$\O_{X}(H)$. In this section we show how $\f(H)$ can be calculated
by a simple analysis of the very same set. The origins of the algorithm
we present here can be traced to \cite{Nic94} with further development
in \cite{LP10}. We present it here from a more general perspective.

Let the group $G$ act on a set $Y$ and let $y_{0}\in Y$ be a fixed
element. Consider a random homomorphism $\alpha_{G}\in Hom(F_{k},G)$.
The core graphs in $\O_{X}(H)$ can be used to calculate the probability
that $\alpha_{G}(H)\subset G_{y_{0}}$ (where $G_{y_{0}}$ is the
stabilizer of the element $y_{0}$). The quotients of the core graph
$\G_{X}(H)$ correspond to all the different {}``coincidence patterns''
of the paths of $y_{0}$ through the action of the images of all $w\in H$,
thereby describing disjoint events whose union is the event that $\a_{G}(H)\subset G_{y_{0}}$.

The idea is that in order to determine whether $\a_{G}(w)$ stabilizes
$y_{0}$ for some $w\in\F_{k}$, we do not need to know all the values
$\a_{G}(x_{i})$ over $x_{i}\in X$ (the given basis of $\F_{k}$).
Rather, we only need to know how $\alpha_{G}(x_{i})$ acts on certain
elements in $Y$, those in the path of $y_{0}$ through $\alpha_{G}(w)$.
Namely, if $w=x_{j_{1}}^{\epsilon_{1}}\ldots x_{j_{|w|}}^{\epsilon_{|w|}}$,
$j_{i}\in\{1,\ldots,k\},\epsilon_{i}\in\{\pm1\}$, we need to uncover
the elements $y_{1},\ldots,y_{|w|}$ in the following diagram:

\[
\xy(0,0)*+{y_{0}}="y0";(20,0)*+{y_{1}}="y1";(40,0)*+{y_{2}}="y2";(60,0)*+{\ldots}="ldots";(80,0)*+{y_{|w|-1}}="y_{n}{}_{1}";(100,0)*+{y_{|w|}}="y_{n}";{\ar^{\a_{G}\big(x_{j_{1}}^{\epsilon_{1}}\big)}"y0";"y1"};{\ar^{\a_{G}\big(x_{j_{2}}^{\epsilon_{2}}\big)}"y1";"y2"};{\ar^{\a_{G}\big(x_{j_{3}}^{\epsilon_{3}}\big)}"y2";"ldots"};{\ar^{\a_{G}\big(x_{j_{|w|-1}}^{\epsilon_{|w|-1}}\big)~~}"ldots";"y_{n}{}_{1}"};{\ar^{\a_{G}\big(x_{j_{|w|}}^{\epsilon_{|w|}}\big)}"y_{n}{}_{1}";"y_{n}"};\endxy
\]
 That is, the image of $x_{j_{1}}^{\epsilon_{1}}$ acts on $y_{0}$,
and we denote the resulting element by $y_{1}\in Y$. The image of
$y_{1}$ under the action of $\a_{G}\big(x_{j_{2}}^{\epsilon_{2}}\big)$
is denoted by $y_{2}$, etc. Then, $y_{0}$ is a fixed point of $\a_{G}(w)$
iff $y_{|w|}=y_{0}$.

There are normally many possible series of elements $y_{1},\ldots,y_{|w|-1}\in Y$
that can constitute the path of $y_{0}$ through $\a_{G}(w)$ such
that $y_{0}$ is a fixed point. We divide these different series to
a \emph{finite} number of categories based on the \emph{coincidence
pattern} of this series. Namely, two realizations of this series,
$y_{1},\ldots,y_{|w|-1}$, and $y'_{1},\ldots,y'_{|w|-1}$ are in
the same category iff for every $i,j\in\{0,\ldots,|w|-1\}$, $y_{i}=y_{j}\Leftrightarrow y'_{i}=y'_{j}$
(note that the elements of the series are also compared to $y_{0}$).
In other words, every coincidence pattern corresponds to some partition
of $\{0,1,\ldots,|w|-1\}$.

However, because the elements $\a_{G}(x_{j})\in G$ act as permutations
on $Y$, not every partition yields a realizable coincidence pattern:
if, for example, $x_{j_{2}}^{\epsilon_{2}}=x_{j_{7}}^{-\epsilon_{7}}$,
and $y_{1}=y_{7}$, we must also have $y_{2}=y_{6}$. This condition
should sound familiar. Indeed, for each coincidence pattern we can
draw a pointed, directed, edge-labeled graph describing it. The vertices
of this graph correspond to blocks in the partition of $\{0,1,\ldots,|w|-1\}$,
the basepoint corresponding to the block containing $0$. Then, for
each $i\in\{1,\ldots,|w|\}$ there is a $j_{i}$-edge, directed according
to $\epsilon_{i}$, between the block of $i-1$ to the block of $i$.
The constraints that coincidence patterns must satisfy then becomes
the very same ones we had encountered in our discussion of core graphs.
Namely, no two $j$-edges share the same origin or the same terminus.

Thus, the different realizable coincidence patterns of the series
$y_{0},y_{1},\ldots,y_{|w|-1}$ are exactly those described by core
graphs that are quotients of $\G_{X}(\la w\ra)$. For instance, there
are exactly seven realizable coincidence patterns that correspond
to the event in which $y_{0}$ is a fixed point of $\a_{G}(w)$ when
$w=[x_{1},x_{2}]$. The seven core graphs in Figure \ref{fig:lattice}
correspond to these seven coincidence patterns.

Finally, the same phenomenon generalizes to any $H\fg\F_{k}$. Instead
of uncovering the path of $y_{0}$ through the image of a single word,
we uncover the paths trough all words in $H$. The union of these
paths in which $y_{0}$ is stabilized by all elements of $H$ is depicted
exactly by the core graph $\G_{X}(H)$. The realizable coincidence
patterns correspond then to the quotients of $\G_{X}(H)$, namely
to $\O_{X}(H)$. To summarize:

\begin{equation}
Prob\big[\a_{G}(H)\subset G_{y_{0}}\big]=\sum_{\G\in\O_{X}(H)}Prob\big[\begin{gathered}\G\textrm{ describes the coincidence pattern }\\
\textrm{of \ensuremath{y_{0}}\,\ through the action of }\a_{G}(H)
\end{gathered}
\big]\label{eq:complement_events}
\end{equation}

The advantage of the symmetric group $S_{n}$ with its action on $\{1,\ldots,n\}$
is that the probabilities in the r.h.s. of \eqref{eq:complement_events}
are very easy to formulate. Let $\a_{n}=\a_{S_{n}}\in Hom(\F_{k},S_{n})$
be a uniformly distributed random homomorphism, and let $\G\in\O_{X}(H)$
be one of the quotients of $\G_{X}(H)$. Denote by $P_{\G}(n)$ the
probability that $\a_{n}(H)\subset(S_{n})_{1}$ and that the coincidence
pattern of the paths of $1$ through the elements $\a_{G}(H)$ are
described by $\G$. Then we can give an exact expression for $P_{\G}(n)$
in terms of $v_{\G},e_{\G}$ and $e_{\G}^{j}$, the number of vertices,
edges and $j$-edges in $\G$:

There are $(n-1)(n-2)\ldots(n-v_{\G}+1)$ possible assignments of
different elements from $\{2,3,\ldots,n\}$ to the vertices of $\G$
(excluding the basepoint which always corresponds to the element $1$).
Then, for a given assignment, there are exactly $e_{\G}^{j}$ constraints
on the permutation $\a_{n}(x_{j})$. So the probability that the permutation
$\a_{n}(x_{j})$ agrees with the given assignment is 
\[
\frac{(n-e_{\G}^{j})!}{n!}=\frac{1}{n(n-1)\ldots(n-e_{\G}^{j}+1)}
\]
 (for $n\ge e_{\G}^{j}$). Thus 
\[
P_{\G}(n)=\frac{(n-1)(n-2)\ldots(n-v_{\G}+1)}{\prod_{j=1}^{k}{n(n-1)\ldots(n-e_{\G}^{j}+1)}}
\]

Recall the definition of $\Phi_{H}(n)$ in \eqref{eq:Phi_def}. Since
for every $j$ and every $\G\in\O_{X}(H)$ we have $e_{\G}^{j}\le e_{\G_{X}(H)}^{j}$
we can summarize and say that for every $n\ge\max_{j}e_{\G_{X}(H)}^{j}$
(in particular for every $n\ge v_{\G_{X}(H)}$), we have:

\begin{align}
\Phi_{H}(n) & =Prob\left[\forall w\in H~~\a_{n}(w)(1)=1\right]-\frac{1}{n^{rk(H)}}\nonumber \\
 & =Prob\left[\a_{n}(H)\subset(S_{n})_{1}\right]-\frac{1}{n^{rk(H)}}\nonumber \\
 & =-\frac{1}{n^{rk(H)}}+\sum_{\G\in\O_{X}(H)}{\frac{(n-1)(n-2)\ldots(n-v_{\G}+1)}{\prod_{j=1}^{k}{n(n-1)\ldots(n-e_{\G}^{j}+1)}}}\nonumber \\
 & =-\frac{1}{n^{rk(H)}}+\sum_{\G\in\O_{X}(H)}\frac{1}{n^{e_{\G}-v_{\G}+1}}\frac{(1-\frac{1}{n})(1-\frac{2}{n})\ldots(1-\frac{v_{\G}-1}{n})}{\prod_{j=1}^{k}{(1-\frac{1}{n})\ldots(1-\frac{e_{\G}^{j}-1}{n})}}\label{eq:Phi_formula}
\end{align}

For instance, for $H=\la[x_{1},x_{2}]\ra$ there are seven summands
in the r.h.s. of \eqref{eq:Phi_formula}, corresponding to the seven
core graphs in Figure \ref{fig:lattice}. If we go over these core
graphs from top to bottom and left to right, we obtain that for every
$n\ge2$: 
\begin{align*}
\Phi_{\la[x_{1},x_{2}]\ra}(n) & =-\frac{1}{n}+\frac{(n-1)(n-2)(n-3)}{n(n-1)\cdot n(n-1)}+\\
 & \qquad{}+\frac{n-1}{n(n-1)\cdot n}+\frac{n-1}{n\cdot n(n-1)}+\frac{(n-1)(n-2)}{n(n-1)\cdot n(n-1)}+\\
 & \qquad{}+\frac{(n-1)(n-2)}{n(n-1)\cdot n(n-1)}+\frac{n-1}{n(n-1)\cdot n(n-1)}+\frac{1}{n\cdot n}\\
 & =-\frac{1}{n}+\frac{1}{n-1}=\frac{1}{n(n-1)}
\end{align*}

Recall the definition of the second categorization of sets of free
words, $\phi(H)$, in \eqref{eq:phi}. Indeed, we can rewrite \eqref{eq:Phi_formula}
as a power series in $\frac{1}{n}$, and obtain that (for large enough
$n$) 
\[
\Phi_{H}(n)=\sum_{i=0}^{\infty}\frac{a_{i}(H)}{n^{i}}
\]
 where the coefficients $a_{i}(H)$ depend only on $H$. We need not
consider negative values of $i$ because the leading term of every
summand in \eqref{eq:Phi_formula} is $\frac{1}{n^{e_{\G}-v_{\G}+1}}$,
and $e_{\G}-v_{\G}+1$ is non-negative for connected graphs. In fact,
this number also equals the rank of the free subgroup corresponding
to $\G$.

The value of $\phi(H)$ equals the smallest $i$ for which $a_{i}(H)$
does not vanish. For instance, for $H=\la[x_{1},x_{2}]\ra$ we have
\[
\Phi_{\la[x_{1},x_{2}]\ra}(n)=\frac{1}{n(n-1)}=\sum_{i=2}^{\infty}\frac{1}{n^{i}}
\]
 so that $a_{0}(H)=a_{1}(H)=0$ and $a_{i}(H)=1$ for $i\geq2$. Hence
$\phi(H)=2$.

In fact, we can write a power series for each $\G\in\O_{X}(H)$ separately,
and obtain: 
\begin{align}
P_{\G}(n) & =\frac{1}{n^{e_{\G}-v_{\G}+1}}\frac{(1-\frac{1}{n})(1-\frac{2}{n})\ldots(1-\frac{v_{\G}-1}{n})}{\prod_{j=1}^{k}{(1-\frac{1}{n})\ldots(1-\frac{e_{\G}^{j}-1}{n})}}\nonumber \\
 & =\frac{1}{n^{e_{\G}-v_{\G}+1}}\big(1-\frac{\binom{v_{\G}}{2}-\sum_{j=1}^{k}\binom{e_{\G}^{j}}{2}}{n}+O\big(\frac{1}{n^{2}}\big)\big)\label{eq:P_G}
\end{align}

This shows that if $\G=\G_{X}(J)$ ($J\fg\F_{k}$), then $P_{\G}(n)$
never affects $a_{i}(H)$-s with $i<rk(J)$. It is also easy to see
that all the coefficients of the power series expressing $P_{\G}(n)$
are integers. We summarize:
\begin{claim}
\label{clm:a_i(S)} For every $H\fg\F_{k}$, all the coefficients
$a_{i}(H)$ are integers. \\
 Moreover, $a_{i}(H)$ is completely determined by core graphs in
$\O_{X}(H)$ corresponding to groups of rank $\le i$. 
\end{claim}

\section{Relations between $\pi(\cdot)$ and $\phi(\cdot)$}

\label{sec:proofs}

We now have all the background needed for the proof of Proposition
\ref{prop:phi=00003Dpi} and consequently of Theorem \ref{ther:m.p.=00003D=00003D>prim._when_r>=00003Dk-1}.
We need to show that for every $H\fg\F_{k}$ and every $i\le rk(H)+1$,
we have 
\[
\pi(H)=i\Longleftrightarrow\f(H)=i.
\]

The proof is divided into three steps. First we deal with the case
$i<rk(H)$, then with $i=rk(H)$. The last case $i=rk(H)+1$ is by
far the hardest.
\begin{lem}
\label{lem:i<rk(H)} Let $H\fg\F_{k}$ and $i<rk(H)$. Then \end{lem}
\begin{enumerate}
\item $\pi(H)=i\Leftrightarrow\f(H)=i$ 
\item If $\pi(H)=\phi(H)=i$ then $a_{i}(H)$ equals the number of $H$-critical
subgroups of $\F_{k}$. \end{enumerate}
\begin{proof}
Let $m$ denote the smallest rank of a group $J\le\F_{k}$ such that
$H\covers J$ (so $m\leq\mathrm{rk\left(H\right)}$). The first part
of the result is derived from the observation that both $\pi(H)=i$
and $\f(H)=i$ iff $m=i$. Let us note first that $\pi(H)=i\Leftrightarrow m=i$.
This follows from Lemma \ref{lem:pi=00003Dpi'} and the fact that
$H$ cannot be a free factor in a subgroup of smaller rank.

We next observe that $\phi(H)=i\Leftrightarrow m=i$: If $m<rk(H)$
then by \eqref{eq:Phi_formula} and \eqref{eq:P_G}, $m$ is indeed
the smallest index for which $a_{m}(H)$ does not vanish (this does
not work for $m=rk(H)$ because of the term $\big(-\frac{1}{n^{rk(H)}}\big)$
in the definition of $\Phi_{H}(n)$). Conversely, if $m=rk(H)$ then
obviously $\phi(H)\ge rk(H)$.

For the second part of the lemma, recall that $H$ is not a free factor
in any subgroup of smaller rank containing it. Thus, by \eqref{eq:P_G}
and Lemma \ref{lem:pi=00003Dpi'}, both $a_{i}(H)$ and the number
of subgroups of rank $i$ containing $H$ equal the number of subgroups
of rank $i$ in $\O_{X}(H)$. 
\end{proof}
The case $i=rk(H)$ is slightly different, but almost as easy.
\begin{lem}
\label{lem:i=00003Drk(H)} Let $H\fg\F_{k}$. Then, \end{lem}
\begin{enumerate}
\item $\pi(H)=rk(H)\Leftrightarrow\f(H)=rk(H)$ 
\item If $\pi(H)=\phi(H)=rk(H)$ then $a_{rk(H)}(H)$ equals the number
of $H$-critical subgroups of $\F_{k}$. \end{enumerate}
\begin{proof}
From Lemma \ref{lem:i<rk(H)} we infer that $\pi(H)\ge rk(H)\Leftrightarrow\f(H)\ge rk(H)$.
So we assume that indeed $\pi(H),\phi(H)\ge rk(H)$, or, equivalently,
that there are no subgroups covered by $H$ of rank smaller than $rk(H)$.

We show that both sides of part $(1)$ are equivalent to the existence
of a quotient (corresponding to a subgroup) of rank $rk(H)$ in $\O_{X}(H)$
\emph{other} than $\G_{X}(H)$ itself. Indeed, this is true for $\pi(H)$
because the only free product of $H$ of rank $rk(H)$ is $H$ itself.

As for $\phi(H)$, this is true because when $\phi(H)\geq rk(H)$
it is easily verified that the value of $a_{rk(H)}(H)$ equals the
number of quotient in $\O_{X}(H)$ of rank $rk(H)$ minus $1$ (this
minus $1$ comes from the term $\big(-\frac{1}{n^{rk(H)}}\big)$).
We think of this term as offsetting the contribution of $\G_{X}(H)$
to $a_{rk(H)}(H)$, so $a_{rk(H)}(H)$ equals the number of other
quotients in $\O_{X}(H)$ of rank $rk(H)$.

The second part of the lemma is true because all $H$-critical subgroups
are covered by $H$ (Lemma \ref{lem:pi=00003Dpi'}).
\end{proof}

\subsection{The Case $i=rk(H)+1$}

\label{sbs:rkH+1} The most interesting (and the hardest) case of
Theorem \ref{ther:m.p.=00003D=00003D>prim._when_r>=00003Dk-1} is
when $rk(H)=k-1$. In the previous analysis this corresponds to $i=rk(H)+1$.
\begin{lem}
\label{lem:i=00003Drk(H)+1} Let $H\fg\F_{k}$. Then, \end{lem}
\begin{enumerate}
\item $\pi(H)=rk(H)+1\Leftrightarrow\f(H)=rk(H)+1$ 
\item If $\pi(H)=\phi(H)=rk(H)+1$ then $a_{rk(H)+1}(H)$ equals the number
of $H$-critical subgroups of $\F_{k}$. 
\end{enumerate}
Denote by $\hG=\G_{X}(H)$ the associated core graph. By Lemmas \ref{lem:i<rk(H)}
and \ref{lem:i=00003Drk(H)}, we can assume that $\pi(H),\phi(H)\geq rk(H)+1$.
In particular, we can thus assume that $H$ is not contained in any
subgroup of rank smaller than $rk(H)+1$ other than $H$ itself.

The coefficient $a_{rk(H)+1}(H)$ in the expression of $\Phi_{H}(n)$
is the sum of two expressions: 
\begin{itemize}
\item The contribution of $\hG$ which equals $-\big(\binom{v_{\hG}}{2}-\sum_{j=1}^{k}\binom{e_{\hG}^{j}}{2}\big)$ 
\item A contribution of $1$ from each core graph of rank $rk(H)+1$ in
$\O_{X}(H)$ 
\end{itemize}
Thus, our goal is to show that the contribution of $\hG$ is exactly
offset by the contribution of the core graphs of rank $rk(H)+1$ in
$\O_{X}(H)$ in which $H$ is a free factor. This would then yield
immediately both parts of Lemma \ref{lem:i=00003Drk(H)+1}. But the
number of subgroups of rank $rk(H)+1$ (in $\O_{X}(H)$) in which
$H$ is a free factor equals exactly the number of immediate quotients
of $\hG$: Theorem \ref{thr:rho=00003Drk-rk_iff_ff} shows that only
immediate quotients of $\hG$ are subgroups of rank $rk(H)+1$ in
which $H$ is a free factor. On the other hand, \eqref{eq:rho>=00003Drk-rk}
and the assumption that $H$ in not contained in any other subgroup
of equal or smaller rank yield that every immediate quotient of $\hG$
is of rank $rk(H)+1$ (and $H$ is a free factor in it).

Thus, Lemma \ref{lem:i=00003Drk(H)+1} follows from the following
lemma.
\begin{lem}
\label{lem:n_immd_quots} Assume $\pi(H),\phi(H)>rk(H)$. Then $\hG=\G_{X}(H)$
has exactly 
\[
\binom{v_{\hG}}{2}-\sum_{j=1}^{k}\binom{e_{\hG}^{j}}{2}
\]
 immediate quotients.
\end{lem}
The intuition behind Lemma \ref{lem:n_immd_quots} is this: Every
immediate quotient of $\hG$ is generated by identifying some pair
of vertices of $\hG$, and there are exactly $\binom{v_{\hG}}{2}$
such pairs. But for every pair of equally-labeled edges of $\hG$,
the pair of origins generates the same immediate quotient as the pair
of termini. This intuition needs, however, some justification that
we now provide.

To this end we use the graph $\U$, a concept introduced in \cite{LP10}%
\footnote{This is a variation of the classical construction of \emph{pull-back}
of graphs (in this case the pull-back of the graph $\hG$ with itself).%
}. This graph represents the pairs of vertices of $\hG$ and the equivalence
relations between them induced by equally-labeled edges. There are
$\binom{v_{\hG}}{2}$ vertices in $\U$, one for each unordered pair
of vertices of $\hG$. The number of directed edges in $\U$ is $\sum_{j=1}^{k}{\binom{e_{\hG}^{j}}{2}}$,
one for each pair of equally-labeled edges in $\hG$. The edge corresponding
to the pair $\{\epsilon_{1},\epsilon_{2}\}$ of $j$-edges is a $j$-edge
connecting the vertex $\{origin(\epsilon_{1}),origin(\epsilon_{2})\}$
to $\{terminus(\epsilon_{1}),terminus(\epsilon_{2})\}$. For example,
when $S$ consists of the commutator word, $\U$ has $\binom{4}{2}=6$
vertices and $\binom{2}{2}+\binom{2}{2}=2$ edges. We illustrate a
slightly more interesting case in Figure \ref{fig:Upsilon}.

\begin{figure}[h]
\begin{centering}
\begin{center}
\xy 
(16,34)*+{\otimes}="v0"+(0,3)*{v_0};
(28,29)*+{\bullet}="v1"+(3,2)*{v_1};
(32,14)*+{\bullet}="v2"+(3,0)*{v_2};
(24,0 )*+{\bullet}="v3"+(2,-3)*{v_3};
(8 ,0 )*+{\bullet}="v4"+(-2,-3)*{v_4};
(0 ,14)*+{\bullet}="v5"+(-3,0)*{v_5}; 
(4 ,29)*+{\bullet}="v6"+(-3,2)*{v_6};
{\ar^{1} "v0";"v1"}; 
{\ar^{1} "v1";"v2"};
{\ar^{2} "v2";"v3"}; 
{\ar^{1} "v3";"v4"}; 
{\ar^{2} "v4";"v5"}; 
{\ar_{1} "v6";"v5"};
{\ar^{2} "v6";"v0"}; 
%
(50,30)*+{\bullet}="t0"+(0,3)*{\scriptstyle \{v_0,v_1\}}; 
(50,15)*+{\bullet}="t1"+(0,-3)*{\scriptstyle \{v_1,v_2\}};
(50,0 )*+{\bullet}="t2"+(0,-3)*{\scriptstyle \{v_3,v_4\}}; 
(64,30)*+{\bullet}="t3"+(0,3)*{\scriptstyle \{v_0,v_2\}}; 
(64,15)*+{\bullet}="t4"+(0,3)*{\scriptstyle \{v_4,v_6\}};
(64,0 )*+{\bullet}="t5"+(0,-3)*{\scriptstyle \{v_0,v_5\}}; 
(78,30)*+{\bullet}="t6"+(0,3)*{\scriptstyle \{v_0,v_4\}}; 
(78,15)*+{\bullet}="t7"+(0,3)*{\scriptstyle \{v_1,v_6\}}; 
(78,0 )*+{\bullet}="t8"+(0,-3)*{\scriptstyle \{v_2,v_5\}}; 
(92,30)*+{\bullet}="t9"+(0,3)*{\scriptstyle \{v_2,v_6\}}; 
(92,15)*+{\bullet}="t10"+(-6,0)*{\scriptstyle \{v_0,v_3\}}; 
(92,0 )*+{\bullet}="t11"+(0,-3)*{\scriptstyle \{v_1,v_4\}}; 
(106,30)*+{\bullet}="t12"+(0,3)*{\scriptstyle \{v_2,v_3\}}; 
(106,15)*+{\bullet}="t13"+(0,3)*{\scriptstyle \{v_3,v_6\}}; 
(106,0 )*+{\bullet}="t14"+(0,-3)*{\scriptstyle \{v_4,v_5\}}; 
(120,30)*+{\bullet}="t15"+(0,3)*{\scriptstyle \{v_0,v_6\}}; 
(120,15)*+{\bullet}="t16"+(0,-3)*{\scriptstyle \{v_1,v_5\}}; 
(120,0 )*+{\bullet}="t17"+(0,-3)*{\scriptstyle \{v_5,v_6\}}; 
(136,30)*+{\bullet}="t18"+(0,3)*{\scriptstyle \{v_1,v_3\}}; 
(136,15)*+{\bullet}="t19"+(6,0)*{\scriptstyle \{v_2,v_4\}}; 
(136,0 )*+{\bullet}="t20"+(0,-3)*{\scriptstyle \{v_3,v_5\}}; 
{\ar_{1} "t0";"t1"}; 
{\ar_{2} "t4";"t5"}; 
{\ar_{1} "t7";"t8"}; 
{\ar_{2} "t9";"t10"};
{\ar_{1} "t10";"t11"};
{\ar_{1} "t13";"t14"};
{\ar_{1} "t15";"t16"}; 
{\ar_{1} "t18";"t19"}; 
{\ar_{2} "t19";"t20"}; 
\endxy 
\par\end{center}
\par\end{centering}

\caption{The graph $\U$ (on the right) corresponding to $\hG=\G_{X}(H)$ (on
the left) for $H=\la x_{1}^{2}x_{2}x_{1}x_{2}x_{1}^{\;-1}x_{2}\ra$.
(The vertices of $\hG$ are denoted here by $v_{0},\ldots,v_{6}$.) }

\label{fig:Upsilon} 
\end{figure}

We denote the set of connected components of $\U$ by $\Co(\U)$.
The proof of Lemma \ref{lem:n_immd_quots} will follow from two facts
that we show next. Namely, $\U$ has exactly $\binom{v_{\hG}}{2}-\sum_{j=1}^{k}\binom{e_{\hG}^{j}}{2}$
connected components. Also, there is a one-to-one correspondence between
$\Co(\U)$ and the set of immediate quotients of $\hG$.
\begin{claim}
\label{clm:n_comps} If $\pi(H),\phi(H)>rk(H)$, then 
\[
\big|\Co(\U)\big|=\binom{v_{\hG}}{2}-\sum_{j=1}^{k}\binom{e_{\hG}^{j}}{2}
\]
\end{claim}
\begin{proof}
Since $\U$ has $\binom{v_{\hG}}{2}$ vertices and $\sum_{j=1}^{k}\binom{e_{\hG}^{j}}{2}$
edges, it is enough to show that it is a forest, i.e., it contains
no cycles.

Let $C\in\Co(\U)$ be some component of $\U$. Clearly, every vertex
in $C$ (which corresponds to a pair of vertices in $\hG$) generates
the same immediate quotient. Denote this quotient by $\Delta(C)$,
and the corresponding subgroup by $J$. In particular, $rk(J)=rk(H)+1$
(recall that under the claim's assumptions, $H$ is not contained
in any other subgroup of smaller or equal rank).Assume to the contrary
that $C$ contains a cycle. Edges in $\U$ are directed and labeled,
so there is an element $u\in\F_{k}$ which corresponds to this cycle
started, say, at the vertex $\{x,y\}$.

\begin{center}
\begin{center}
\xy 
(40,30)*+{\bullet}="v0"+(0,4)*{\{x,y\}}; 
"v0";"v0" **\crv{~*=<4pt>{.} (60,20) & (40,0) & (20,20) } ?<(1)*\dir{>}; 
(54,20)*+{u};
(-30,10)*+{}="dummy";
(-30,33)*+{}="dummy";
\endxy
\par\end{center}
\par\end{center}

Where do we get as we walk in the core graph $\hG$ starting at $x$
(resp. $y$) and following the path corresponding to $u$? One possibility
is that the walk from $x$ returns back to $x$ and likewise for $y$.
Alternatively this $u$-walk can take us from $x$ to $y$ and from
$y$ to $x$. We consider only the former possibility. The latter
case would be handled by considering the walk corresponding to $u^{2}$.
Let $p_{x},p_{y}\in\F_{k}$ be words corresponding to some paths from
the basepoint of $\hG$ to $x,y$ respectively. In particular, $p_{x}up_{x}^{-1},p_{y}up_{y}^{-1}\in H$.

\begin{center}
\begin{center}
\xy
(0,0)*{\otimes}="bp"+(-8,8)*{\hG};%
(5,0)*\xycircle(30,15){--};
(24,4)*{\bullet}="x"  +(-1,3)*{x};%
(20,-6)*{\bullet}="y" +(-2,-2)*{y};%
{\ar@(dr,ur)_{u} "x";"x"};%
{\ar@(dr,ur)_{u} "y";"y"};%
"bp"; "x" **\crv{(10,20) & (15,-20)} ?(.43)*\dir{>>}; 
(6,8)*{\scriptstyle p_x}; 
"bp"; "y" **\crv{(5,-10) & (15,10)} ?(.43)*\dir{>>}; 
(6,-6)*{\scriptstyle p_y}; 
(-60,0)*+{}="dummy";
(0,18)*+{}="dummy";
(0,-18)*+{}="dummy";
\endxy 
\par\end{center}
\par\end{center}

Merging $x$ and $y$ is equivalent to adding the generator $p_{x}p_{y}^{-1}$
to $H$, so that $J=\la H,p_{x}p_{y}^{-1}\ra$. Since $rk(J)=rk(H)+1$,
we have that $J=H*\la p_{x}p_{y}^{-1}\ra$. Consider the elements
$h_{1}=p_{x}up_{x}^{-1}\in H$ and $h_{2}=p_{y}up_{y}^{-1}\in H$.
The following equality holds: 
\[
h_{1}=p_{x}up_{x}^{-1}=(p_{x}p_{y}^{-1})p_{y}up_{y}^{-1}(p_{x}p_{y}^{-1})^{-1}=(p_{x}p_{y}^{-1})h_{2}(p_{x}p_{y}^{-1})^{-1}
\]
 This is a contradiction, since we obtained two different expressions
for $h_{1}$ in the free product $J=H*\la p_{x}p_{y}^{-1}\ra$.\end{proof}
\begin{rem}
Let $H\leq\F_{k}$ and $x,y\in V\left(\Gamma_{X}\left(H\right)\right)$.
In the proof of the last claim it was shown that if there is some
$1\ne u\in\F_{k}$ which is readable as a closed path at both $x$
and $y$, then the subgroup we obtain by merging them is \emph{not}
a free extension of $H$. We stress that the converse is not true.
For example, consider $H=\left\langle a,bab\right\rangle \leq\F\left(\left\{ a,b\right\} \right)$.
Then $\Gamma_{\left\{ a,b\right\} }\left(H\right)$ has three vertices,
no pair of which share a common closed path (in other words, the corresponding
graph $\Upsilon$ has no cycles). However, $\Gamma_{\left\{ a,b\right\} }\left(H\right)$
has exactly two immediate quotients, none of which is a free extension.
\end{rem}
Next we exhibit a one-to-one correspondence between $\Co(\U)$ and
the immediate quotients of $\hG$. It is very suggestive to try and
restore $C$ from $\Delta(C)$ by simply signaling out the pairs of
vertices that are identified in $\Delta(C)$. But this is too naive.
There may be pairs of vertices not in $C$ that are identified in
$\Delta(C)$. For instance, consider $C$, the rightmost component
of $\U$ in Figure \ref{fig:Upsilon}. In $\Delta(C)$ we merge $v_{1}$
and $v_{3}$ but also $v_{3}$ and $v_{5}$. Thus $v_{1}$ and $v_{5}$
are merged and likewise all pairs in the component of $\{v_{1},v_{5}\}$.

However, simple group-theoretic arguments do yield this sought-after
result: 
\begin{claim}
\label{clm:correspondence} If $\pi(H),\phi(H)>rk(H)$, then there
is a one-to-one correspondence between $\Co(\U)$ and the set of immediate
quotients of $\hG=\G_{X}(H)$.\end{claim}
\begin{proof}
Maintaining the above notation, we need to show that the mapping from
$C\in\Co(\U)$ to $\Delta(C)$, the immediate quotient generated by
any of the pairs in $C$, is one to one.

Denote by $J$ the subgroup corresponding to the immediate quotient
$\Delta(C)$. Let $\{x,y\}$ be some vertex in $C$, and $p_{x},p_{y}\in\F_{k}$
words corresponding to some paths from the basepoint of $\hG$ to
$x,y$, respectively. Let also $q=p_{x}p_{y}^{-1}\in\F_{k}$. As we
saw above, 
\[
J=\langle H,q\rangle
\]
 and clearly $q\notin H$.

We claim that any other complementary generator of $J$ over $H$
is in same $(H,H)$-double-coset of $q$ or of $q^{-1}$ in $J$.
Namely, if $J=\langle H,q'\rangle$ then $q'\in HqH\cup Hq^{-1}H$.
To see this, let $Y$ be some basis of $H$ and think of $J$ as the
free group over the basis $Y\cup\{q\}$ (this is true because $rk(J)=rk(H)+1$).
Now think of $q'$ as a word in the elements of this basis. Multiplying
from the right or left by elements of $Y$ does not affect the $(H,H)$-double-coset,
so assume w.l.o.g. that $q'$ begins and ends with either $q$ or
$q^{-1}$. But then the set $Y\cup\{q'\}$ is Nielsen-reduced with
respect to the basis $Y\cup\{q\}$ (see, for instance, the definition
in Chapter 1 of \cite{LS70}). As consequence, $Y\cup\{q'\}$ equals
$Y\cup\{q\}$ up to taking inverses (Proposition 2.8 therein). Thus
$q'=q$ or $q'=q^{-1}$.

So let $\{a,b\}$ be another pair of vertices generating $\Delta(C)$.
We show that it belongs to $C$. Let $p_{a},p_{b}$ be words in $\F_{k}$
corresponding to paths in $\hG$ from the basepoint to $a,b$ respectively.
We have $\langle H,p_{a}p_{b}^{-1}\rangle=J$, so $p_{a}p_{b}^{-1}\in HqH\cup Hq^{-1}H$.
W.l.o.g. it is in $HqH$ (otherwise switch $a$ and $b$). So assume
$p_{a}p_{b}^{-1}=h_{1}qh_{2}$ with $h_{1},h_{2}\in H$. But $h_{1}^{-1}p_{a}$
is also a path from the basepoint of $\hG$ to $a$, and likewise
$h_{2}p_{b}$ a path to $b$. Choosing if needed these paths instead
of $p_{a},p_{b}$ we can assume that 
\[
p_{a}p_{b}^{-1}=q=p_{x}p_{y}^{-1}.
\]

Thus, 
\[
p_{a}^{-1}p_{x}=p_{b}^{-1}p_{y}.
\]
 This shows that there is a path in $\hG$ from $a$ to $x$ corresponding
to a path from $b$ to $y$. This shows precisely that the pair $\{a,b\}$
is in the same component of $\U$ as $\{x,y\}$, namely, in $C$. 
\end{proof}
This completes the proof of Lemma \ref{lem:i=00003Drk(H)+1}. This
Lemma, together with Lemmas \ref{lem:i=00003Drk(H)} and \ref{lem:i<rk(H)},
yields Proposition \ref{prop:phi=00003Dpi} and thus Theorem \ref{ther:m.p.=00003D=00003D>prim._when_r>=00003Dk-1}.

\subsection{Further Relations between $\pi(\cdot)$ and $\phi(\cdot)$}

\label{sbs:more_relations}

Let us take another look now at Conjecture \ref{conj:phi=00003Dpi}.
It posits that the results described in Proposition \ref{prop:phi=00003Dpi}
hold for all values of $\pi(\cdot)$ and $\phi(\cdot)$. To understand
what this means, suppose that $H$ is a free factor in all the quotients
in $\O_{X}(H)$ of ranks up to $i-1$. What can be said about rank-$i$
quotients in which $H$ is a free factor? The conjecture states that
their number exactly offsets the sum of two terms: The contribution
to $a_{i}(H)$ of the quotients of smaller rank and of the term $\frac{-1}{n^{rk(H)}}$
when $i=rk(H)$. For instance, $\pi(H)=3$ for $H=\la x_{1}^{\;2}x_{2}^{\;2}x_{3}^{\;2}\ra$.
In particular, $H$ is a free factor of all quotients in $\O_{X}(H)$
of rank $\le2$. There is a single $H$-critical subgroup ($\F_{3}$
itself), and additional $13$ quotients of rank $3$ in which $H$
is a free factor. The contribution of quotients of rank $\le2$ to
$a_{3}(H)$ is indeed exactly $(-13)$.

Interestingly enough, this is indeed the case for every free factor
$H\ff\F_{k}$. In this case, since free factors are measure preserving,
we get that $\phi(H)=\infty$, so $a_{i}(H)=0$ for every $i$, and
the statement of the previous paragraph holds. For the general case
the conjecture states that as long as we consider low-rank quotients
and {}``imprimitivity has not been revealed yet'', the situation
does not differ from what is seen in the primitive case. \\

We finish this section by stating another result connecting $\pi(\cdot)$
and $\phi(\cdot)$. It shows an elegant property of both of them that
lends further support to our belief in Conjecture \ref{conj:phi=00003Dpi}.

Two words $w_{1},w_{2}\in\F_{k}$ are called \emph{disjoint} (with
respect to a given basis) if they share no common letters.
\begin{lem}
\label{lem:additivity} Let $w_{1},w_{2}\in\F_{k}$ be disjoint. Then
\begin{eqnarray*}
\pi(w_{1}w_{2})=\pi(w_{1})+\pi(w_{2})\\
\phi(w_{1}w_{2})=\phi(w_{1})+\phi(w_{2})
\end{eqnarray*}
 Moreover, $a_{\phi(w_{1}w_{2})}(w_{1}w_{2})=a_{\phi(w_{1})}(w_{1})\cdot a_{\phi(w_{2})}(w_{2})$,
and if part 2 of Conjecture \ref{conj:phi=00003Dpi} holds for $H=\la w_{1}\ra$
and for $H=\la w_{2}\ra$, then it also holds for $H=\la w_{1}w_{2}\ra$. 
\end{lem}
This lemma is essentially outside the scope of the present paper,
so we only sketch its proof. Let $\alpha_{n}\in Hom(\F_{k},S_{n})$
be a random homomorphism chosen with uniform distribution. As $w_{1}$
and $w_{2}$ are disjoint, the random permutations $\alpha_{n}(w_{1})$
and $\alpha_{n}(w_{2})$ are independent. The claims about the additivity
of $\phi(\cdot)$ and the multiplicativity of $a_{\phi(\cdot)}(\cdot)$
are easy to derive by calculating the probability that $1$ is a fixed
point of $w_{1}w_{2}$. The key fact in this calculation is the aforementioned
independence of $\alpha_{n}(w_{1})$ and $\alpha_{n}(w_{2})$.

The other claims in the lemma follow from an analysis of $H$-critical
subgroups. By considering properties of the associated core graphs
it is not hard to show that $J\le\F_{k}$ is $\la w_{1}w_{2}\ra$-critical
iff it is the free product of a $\la w_{1}\ra$-critical subgroup
and a $\la w_{2}\ra$-critical subgroup.

\section{Primitive Words and the Profinite Completion}

\label{sec:profinite} Most of the standard facts below about profinite
groups and particularly free profinite groups can be found with proofs
in \cite{Wil98} (in particular Section 5.1).

A profinite group is a topological group $G$ with any of the following
equivalent properties: 
\begin{itemize}
\item $G$ is the inverse limit of an inverse system of finite groups. 
\item $G$ is compact, Hausdorff and totally disconnected. 
\item $G$ is isomorphic (as a topological group) to a closed subgroup of
a Cartesian product of finite groups. 
\item $G$ is compact and $\bigcap(N|N\triangleleft_{O}G)=1$ 
\end{itemize}
The \emph{free profinite group} on a finite set $X$ is a profinite
group $F$ together with a map $j:X\to F$ with the following universal
property: whenever $\xi:X\to G$ is a map to a profinite group $G$,
there is a unique (continuous) homomorphism $\bar{\xi}:F\to G$ such
that $\xi=\bar{\xi}j$. Such $F$ exists for every $X$ and is unique
up to a (continuous) isomorphism. We call $j(X)$ a basis of $F$.
It turns out that every two bases of $F$ have the same size which
is called the \emph{rank} of $F$. The free profinite group of rank
$k$ is denoted by $\hat{\F}_{k}$. An element $w\in\hF_{k}$ is \emph{primitive}
if it belongs to some basis.

It is a standard fact that $\hF_{k}$ is the profinite completion
of $\F_{k}$ and $\F_{k}$ is naturally embedded in $\hF_{k}$. Moreover,
every basis of $\F_{k}$ is then also a basis for $\hF_{k}$, so a
primitive word $w\in\F_{k}$ is also primitive as an element of $\hF_{k}$.
It is conjectured that the converse also holds:
\begin{conjecture}
\label{conj:profinite} A word $w\in\F_{k}$ is primitive in $\hF_{k}$
iff it is primitive in $\F_{k}$.
\end{conjecture}
This conjecture, if true, immediately implies the following one:
\begin{conjecture}
\label{conj:closed-profinite}The set of primitive elements in $\F_{k}$
form a closed set in the profinite topology.
\end{conjecture}
Conjecture \ref{conj:prim<=00003D>m.p} implies these last two conjectures
(it is in fact equivalent to Conjecture \ref{conj:profinite}, see
below): we define measure preserving elements in $\hF_{k}$ as before.
Namely, an element $w\in\hF_{k}$ is measure preserving if for any
finite group $G$ and a uniformly distributed random (continuous)
homomorphism $\hat{\alpha}_{G}\in Hom(\hF_{k},G)$, the image $\hat{\alpha}_{G}(w)$
is uniformly distributed in $G$. Clearly, an element of $\F_{k}$
is measure preserving w.r.t $\F_{k}$ iff this holds w.r.t. $\hF_{k}$.

As in the abstract case, a primitive element of $\hF_{k}$ is measure
preserving. Conjecture \ref{conj:prim<=00003D>m.p} would therefore
imply that if $w\in\F_{k}$ is primitive in $\hF_{k}$, then $w$
is also primitive w.r.t. $\F_{k}$. In particular, Theorem \ref{ther:m.p.=00003D=00003D>prim._when_r>=00003Dk-1}
yields:
\begin{cor}
\label{cor:prim_in_hat_F_k} Let $S\subset\F_{k}$ be a finite subset
of cardinality $|S|\ge k-1$. Then, 
\[
S\textrm{ can be extended to a basis in }\hF_{k}\Longleftrightarrow S\textrm{ can be extended to a basis in }\F_{k}
\]
 In particular, for every $w\in\F_{2}$: 
\[
w\textrm{ is primitive in }\hF_{2}\Longleftrightarrow w\textrm{ is primitive in }\F_{2}
\]

\end{cor}
This corollary yields, in turn, Corollary \ref{cor:prim_is_closed},
which states the special case of Conjecture \ref{conj:closed-profinite}
for $\F_{2}$.

As shown by Chen Meiri (unpublished), Conjectures \ref{conj:profinite}
and \ref{conj:prim<=00003D>m.p} are equivalent. With his kind permission
we explain this result in this section. Meiri showed that in $\hF_{k}$
primitivity and measure preservation are equivalent (Proposition \ref{prop:chen}
below). Thus, $w\in\F_{k}$ is primitive as an element of $\hF_{k}$
iff it is measure preserving.
\begin{prop}
\label{prop:chen} {[}C. Meiri, unpublished{]} Let $w$ belong to
$\hF_{k}$. Then 
\[
w\textrm{ is primitive }\Longleftrightarrow w\textrm{ is measure preserving }
\]
 \end{prop}
\begin{proof}
The $(\Rightarrow)$ implication is trivial as in the abstract case:
for every finite group $G$ and every basis $x_{1},\ldots,x_{k}$
of $\hF_{k}$ there is a bijection 
\begin{eqnarray*}
Hom(\hF_{k},G)\stackrel{\cong}{\to}G^{k}\\
\alpha_{G}\mapsto\left(\alpha_{G}(x_{1}),\ldots,\alpha_{G}(x_{k})\right)
\end{eqnarray*}
For the other direction, for every $w\in\hF_{k}$, finite group $G$
and $g\in G$ define 
\begin{eqnarray*}
H_{w}(G,g)=\left\{ \alpha_{G}\in Hom(\hF_{k},G)~|~\alpha_{G}(w)=g\right\} \\
E_{w}(G,g)=\left\{ \alpha_{G}\in Epi(\hF_{k},G)~|~\alpha_{G}(w)=g\right\} 
\end{eqnarray*}
Now assume $w\in\hF_{k}$ is measure preserving, and let $x\in\hF_{k}$
be any primitive element. For every finite group $G$ we have $|H_{w}(G,g)|=|G|^{k-1}=|H_{x}(G,g)|$.
The same equality holds for the set of epimorphisms, namely $|E_{w}(G,g)|=|E_{x}(G,g)|$.
We will show this by induction on $|G|$.

If $|G|=1$ the claim is trivial. The inductive step goes as follows:
if $g\in G$, then 
\begin{eqnarray*}
|E_{w}(G,g)| & = & |H_{w}(G,g)|-\sum_{g\in H\lneqq G}|E_{w}(H,g)|=\\
 & = & |H_{x}(G,g)|-\sum_{g\in H\lneqq G}|E_{x}(H,g)|=|E_{x}(G,g)|
\end{eqnarray*}
Now choose a basis $x_{1},\ldots,x_{k}$ of $\hF_{k}$. For every
$N\triangleleft_{O}\hF_{k}$, $|E_{x_{1}}(\hF_{k}/N,wN)|=|E_{w}(\hF_{k}/N,wN)|\ge1$.
If $\alpha\in E_{x_{1}}(\hF_{k}/N,wN)$ then $wN=\alpha(x_{1}),\alpha(x_{2}),\ldots,\alpha(x_{k})$
generate $\hF_{k}/N$. A standard compactness argument shows that
there are elements $w_{2},\ldots,w_{k}\in\hF_{k}$ such that $\{wN,w_{2}N,\ldots,w_{k}N\}$
generate $\hF_{k}/N$ for every $N\triangleleft_{O}\hF_{k}$. But
then $\{w,w_{2},\ldots,w_{k}\}$ generate $\hF_{k}$ as well. Whenever
$k$ elements generate $\hF_{k}$, they generate it freely. Thus $\{w,w_{2},\ldots,w_{k}\}$
is a basis and $w$ is primitive. 
\end{proof}

\section{The Average Number of Fixed Points in $\alpha_{n}(w)$}

\label{sec:fixed_points} As before, let $\alpha_{n}\in Hom(\F_{k},S_{n})$
be a uniformly distributed random homomorphism. In \eqref{eq:Phi_def}
we defined the function $\Phi_{\la w\ra}(n)=\Phi_{w}(n)$ for every
$w\in\F_{k}$. It considers the probability that $\alpha_{n}(w)$
fixes the element $1$ and quantifies its deviation from $\frac{1}{n}$.
The choice of the element $1$ is arbitrary, of course, and we get
the same probability for every element in $1,\ldots,n$. Thus $n\Phi_{w}(n)+1$
is the average number of fixed points of the random permutation $\alpha_{n}(w)$.

Corollary \ref{cor:image-of-pi} states that in $\F_{2}$ there are
exactly four possible primitivity ranks of words. This translates
through Proposition \ref{prop:phi=00003Dpi} to four possibilities
for the average number of fixed points in the permutation $\alpha_{n}(w)$,
as summarized by Table \ref{tab:fixed-points-F2}:

\begin{table}[htb]
 \centering %
\begin{tabular}{|c|c|c|c|}
\hline 
$\pi(w)/\phi(w)$  & Description  & $Prob[\alpha_{n}(w)(1)=1]$  & Avg \# of f.p. of $\alpha_{n}(w)$ \tabularnewline
\hline 
\hline 
0  & $w=1$  & $1$  & $n$ \tabularnewline
\hline 
1  & $w$ is a power  & $\frac{1}{n}+\frac{a_{1}(w)}{n}+\sum_{i=2}^{\infty}\frac{a_{i}(w)}{n^{i}}$  & $1+a_{1}(w)+O\big(\frac{1}{n}\big)$ \tabularnewline
\hline 
2  &  & $\frac{1}{n}+\frac{a_{2}(w)}{n^{2}}+\sum_{i=3}^{\infty}\frac{a_{i}(w)}{n^{i}}$  & $1+\frac{a_{2}(w)}{n}+O\big(\frac{1}{n^{2}}\big)$ \tabularnewline
\hline 
$\infty$  & $w$ is primitive  & $\frac{1}{n}$  & 1 \tabularnewline
\hline 
\end{tabular}\caption{The possibilities for the average number of fixed points of the permutation
$\alpha_{n}(w)$ for some $w\in\F_{2}$.}

\label{tab:fixed-points-F2} 
\end{table}

Recall that all coefficients $a_{i}(w)$ are integers (Claim \ref{clm:a_i(S)}).
Moreover, in these cases $a_{\phi(w)}(w)$ counts the $\la w\ra$-critical
subgroups of $\F_{2}$, so in particular $a_{\phi(w)}(w)>0$. We thus
obtain
\begin{cor}
\label{cor:fixed-points-F_2} For every word $w\in\F_{2}$ and every
large enough $n$, the average number of fixed points of $\alpha_{n}(w)$
is at least $1$. 
\end{cor}
This leads to the following conjecture, which is a consequence of
Conjecture \ref{conj:phi=00003Dpi}:
\begin{conjecture}
\label{conj:fixed-points} For every word $w\in\F_{k}$ and every
large enough $n$, the average number of fixed points of $\alpha_{n}(w)$
is at least $1$. 
\end{conjecture}
Proposition \ref{prop:phi=00003Dpi} says something about free words
in general. If $\phi(w)\le2$ for some $w\in\F_{k}$, then the first
non-vanishing coefficient $a_{\phi(w)}(w)$ is positive. Thus,
\begin{cor}
\label{cor:fixed-points-F_k} For every word $w\in\F_{k}$ the average
number of fixed points in $\alpha_{n}(w)$ is at least $1-O\big(\frac{1}{n^{2}}\big)$. 
\end{cor}
It is suggestive to ask whether Conjecture \ref{conj:fixed-points}
holds for \emph{all} $n$. Namely, is it true that for every $w\in\F_{k}$
and every $n$, the average number of fixed points in $\alpha_{n}(w)$
is at least $1$? By results of Abért (\cite{Abe06}), this statement
turns out to be incorrect.

\section*{A Note Added in Proof}
\begin{rem}
\label{rem:silva-weil1} After this paper was completed, we learned
about the algorithm of Silva and Weil to detect free-factor subgroups
in the free group \cite{SW08}. In essence, their algorithm relies
on the same phenomenon that we independently noticed here. However,
our reasoning is very different, and offers several substantial advantages
over the presentation in \cite{SW08}. A more elaborate discussion
of the differences between the two approaches appears in Appendix
\ref{sec:algo}. 
\end{rem}

\begin{rem}
In subsequent joint work with O. Parzanchevski \cite{PP12}, we manage
to prove Conjecture \ref{conj:prim<=00003D>m.p} in full. That proof
relies on Theorem \ref{thr:rho=00003Drk-rk_iff_ff} and follows the
general strategy laid out in the current paper. In particular, we
establish Conjectures \ref{conj:phi=00003Dpi}, \ref{conj:profinite},
\ref{conj:closed-profinite} and \ref{conj:fixed-points}.
\end{rem}

\section*{Acknowledgements}

It is a pleasure to thank Nati Linial for his support, encouragement
and useful comments. We are also grateful to Aner Shalev for supporting
this research and for his valuable suggestions. We would also like
to thank Tsachik Gelander, Michael Larsen, Alex Lubotzky, Chen Meiri,
Ori Parzanchevski, Iddo Samet and Enric Ventura for their beneficial
comments. We would also like to express our gratefulness to the anonymous
referee for his many valuable comments.

\section*{Appendices}

\begin{appendices}

\section{An Algorithm to Detect Free Factors\label{sec:algo}}

One of the interesting usages of Theorem \ref{thr:rho=00003Drk-rk_iff_ff}
is an algorithm to detect free factor subgroups and consequently,
also to detect primitive words in $\F_{k}$. The algorithm receives
as input $H$ and $J$, two finitely generated subgroups of $\F_{k}$,
and determines whether $H\ff J$. The subgroups $H$ and $J$ are
given to us by specifying a generating set, where members of the generating
sets are words in the elements of the fixed basis $X$. (Note that
the algorithm in particular decides as well whether $H\le J$, but
this is neither hard nor new).

We should mention that ours is not the first algorithm, nor the first
graph-theoretic one, for this problem (see Chapter I.2 in \cite{LS70}).
We already mentioned (Remark \ref{rem:silva-weil1}) \cite{SW08},
who noticed the basic phenomenon underlying our algorithm, albeit
in a very different language. See Remark \ref{rem:silva-weil2} below
for an explanation of the differences. A well-known algorithm due
to Whitehead solves a much more general problem. Namely, for given
$2r$ words $w_{1},\ldots,w_{r},u_{1},\ldots,u_{r}\in\F_{k}$, it
determines whether there is an automorphism $\alpha\in Aut(\F_{k})$
such that $\alpha(w_{i})=u_{i}$ for each $i$ (\cite{Whi36a},\cite{Whi36b}.
For a good survey see Chapter I.4 at \cite{LS70}. A nice presentation
of the restriction of Whitehead's algorithm to our problem appears
in \cite{Sta99}). Quite recently, Roig, Ventura and Weil introduced
a more clever version of the Whitehead algorithm for the case of detecting
primitive words and free factor subgroups \cite{RVW07}. Their version
of the algorithm has polynomial time in both the length of the given
word $w$ (or the total length of generators of a given subgroup $H$)
and in $k$, the rank of the ambient group $\F_{k}$. To the best
of our knowledge, their algorithm is currently the best one for this
problem, complexity-wise. The algorithm we present is, at least naively,
exponential, as we show below (Remark \ref{rem:complexity}). \\

So assume we are given two subgroups of finite rank of $\F_{k}$,
$H$ and $J$, by means of finite generating sets $S_{H},S_{J}$.
Each element of $S_{H},S_{J}$ is assumed to be a word in the letters
$X\cup X^{-1}$ (recall that $X=\{x_{1},\ldots,x_{k}\}$ is the given
basis of $\F_{k}$). To find out whether $H\ff J$, follow the following
steps.\\

\noindent \textbf{Step 1: Construct Core Graphs and Morphism}

\noindent First, construct the core graphs $\G=\G_{X}(H)$ and $\Delta=\G_{X}(J)$
by the process described in Appendix \ref{sec:Stallings-Folding-Algorithm}.
Then, seek a morphism $\eta:\G\to\Delta$. This is a simple process
that can be done inductively as follows: $\eta$ must map the basepoint
of $\G$ to the basepoint of $\Delta$. Now, as long as $\eta$ is
not fully defined, there is some $j$-edge $e=(u,v)$ in $E(\G)$
for which the image is not known yet, but the image of one of the
end points, say $\eta(u)$, is known (recall that $\G$ is connected).
There is at most one possible value that $\eta(e)$ can take, since
the star of $\eta(u)$ contains at most one outgoing $j$-edge. If
there is no such edge, we get stuck. Likewise, $\eta(v)$ must equal
the terminus of $\eta(e)$, and if $\eta(v)$ was already determined
in an inconsistent way, we get a contradiction. If in this process
we never get stuck and never reach a contradiction, then $\eta$ is
defined. If this process cannot be carried out, then there is no morphism
from $\G$ to $\Delta$, and hence $H$ is not a subgroup of $J$
(see Claim \ref{clm:morphism-properties}).\\

\noindent \textbf{Step 2: Reduce to Two Groups with $H\covers J'$}

\noindent After constructing the morphism $\eta:\G\to\Delta$, we
obtain a new graph from $\Delta$ by omitting all edges and all vertices
not in the image of $\eta$. Namely, 
\[
\Delta':=\eta(\G)
\]
 It is easy to see that $\Delta'$ is a core-graph, and we denote
by $J'$ the subgroup corresponding to $\Delta'$. Obviously, $\Delta'$
is a quotient of $\G$, so $H\covers J'$. Moreover, it follows from
Claim \ref{clm:free_factors} that 
\[
H\ff J\Longleftrightarrow H\ff J'.
\]

\noindent \textbf{Step 3: Use $\rho_{X}(H,J')$ to determine whether
$H\ff J'$}

\noindent Now calculate $\rho_{X}(H,J')$ (this is clearly doable
because the subgraph of $\D_{k}$ consisting of quotients of $\G$
is finite). Thanks to Theorem \ref{thr:rho=00003Drk-rk_iff_ff}, $\rho_{X}(H,J')$
determines whether or not $H\ff J'$, and consequently, whether or
not $H\ff J$. \\

\begin{rem}
\label{rem:complexity}The complexity of this algorithm is roughly
$O(v^{2t})$, where $v$ is the number of vertices in $\G_{X}(H)$
and $t$ is the difference in ranks: $t=rk(J)-rk(H)$. Naively, we
need to go over roughly all possible sets of $t$ pairs of vertices
of $\G_{X}(H)$ and try to merge them (see Remark \prettyref{remark:pairs-of-vertices}).
The number of possibilities is at most $\binom{\binom{v}{2}}{t}$,
which shows the claimed bound. (In fact, we can restrict to pairs
where both vertices are in the same fiber of the morphism $\eta:\G_{X}(H)\to\G_{X}(J)$.)
\end{rem}

\subsection{Examples\label{sbs:examples}}

We illustrate the different phases of the algorithm by two concrete
examples. Consider first the groups $H=\langle x_{1}x_{2}x_{1}^{\;-1}x_{2}^{\;-1},x_{2}x_{1}^{\;2}\rangle$
and $J=\langle x_{1}^{\;3},x_{2}^{\;3},x_{1}x_{2}^{-1},x_{1}x_{2}x_{1}\rangle$,
both in $\F_{2}$. The core graphs of these groups are:

\begin{center}
\begin{center}
\xy
( 0,25)*+{\otimes}="s0";%
(25,25)*+{\bullet}    ="s1";%
(25, 0)*+{\bullet}    ="s2";%
( 0, 0)*+{\bullet}    ="s3";%
{\ar^{1} "s0";"s1"};%
{\ar^{2} "s1";"s2"};%
{\ar^{1} "s3";"s2"};%
{\ar^{2} "s0";"s3"};%
{\ar^{1} "s2";"s0"};%
(65,15)*+{\otimes}="r0";
(85,15)*+{\bullet}="r1";
(105,15)*+{\bullet}="r2";
{\ar@/^1pc/_{1} "r0";"r1"};
{\ar@/^1pc/_{1} "r1";"r2"};
{\ar@/_1pc/^{2} "r1";"r2"};
{\ar@/_1pc/^{2} "r0";"r1"};
{\ar@/_3pc/^{1} "r2";"r0"};
{\ar@/^3pc/_{2} "r2";"r0"};
(-20,20)*+{}="dummy";
(0,30)*+{}="dummy";
(0,-5)*+{}="dummy";
\endxy 
\par\end{center}
\par\end{center}

In this case, a morphism $\eta$ from $\G=\G_{X}(H)$ to $\Delta=\G_{X}(J)$
can be constructed. All the vertices of $\G$ are in the image of
$\eta$, and only one edge, the long $2$-edge at the bottom, is not
in $\eta(E(\Gamma))$. Thus $\Delta'$ is:

\begin{center}
\begin{center}
\xy 
(0,15)*+{\otimes}="r0"; 
(20,15)*+{\bullet}="r1"; 
(40,15)*+{\bullet}="r2";
{\ar@/^1pc/_{1} "r0";"r1"}; 
{\ar@/^1pc/_{1} "r1";"r2"}; 
{\ar@/_1pc/^{2} "r1";"r2"};
{\ar@/_1pc/^{2} "r0";"r1"}; 
{\ar@/_3pc/^{1} "r2";"r0"};
(-40,15)*+{}="dummy";
(0,30)*+{}="dummy";
(0,10)*+{}="dummy";
\endxy 
\par\end{center}
\par\end{center}

\noindent and $J'$ is the corresponding subgroup $J'=\la x_{1}^{\;3},x_{1}x_{2}^{-1},x_{1}x_{2}x_{1}\ra$.

Finally, $rk(H)=1-\chi(\G)=2$ and $rk(J')=1-\chi(\Delta')=3$, and
so the difference is $rk(J')-rk(H)=1$. It can be easily verified
that $\Delta'$ is indeed an immediate quotient of $\G$: simply merge
the upper-right vertex of $\G$ with the bottom-left one to obtain
$\Delta'$. Thus $\rho_{X}(H,J')=1=rk(J')-rk(H)$, and so $H\ff J'$
hence $H\ff J$.\\

As a second example, consider the commutator word $w=x_{1}x_{2}x_{1}^{-1}x_{2}^{-1}$.
We want to determine whether it is primitive in $\F_{3}$. We take
$H=\langle w\rangle$ and the core graphs are then

\begin{center}
\begin{center}
\xy 
( 0,25)*+{\otimes}="s0";%
(25,25)*+{\bullet}="s1";%
(25, 0)*+{\bullet}="s2";%
( 0, 0)*+{\bullet}="s3";%
{\ar^{1} "s0";"s1"};%
{\ar^{2} "s1";"s2"};%
{\ar^{1} "s3";"s2"};%
{\ar^{2} "s0";"s3"};%
(85,15)*+{\otimes}="r";
{\ar@(l,dl)_{1} "r";"r"};
{\ar@(r,dr)^{2} "r";"r"}; 
{\ar@(ul,ur)^{3} "r";"r"};
(-20,20)*+{}="dummy";
(0,29)*+{}="dummy";
(0,-4)*+{}="dummy";
\endxy 
\par\end{center}
\par\end{center}

Once again, a morphism $\eta$ from $\G=\G_{X}(H)$ to $\Delta=\G_{X}(\F_{3})$
can be constructed, and there is a single edge in $\Delta$, the $3$-edge,
outside the image of $\eta$. Thus $\Delta'$ is the quotient of $\G$
which is the bottom graph in Figure \ref{fig:lattice}, and $J'$
is simply $\F_{2}$.

Finally, $rk(H)=1-\chi(\G)=1$ and $rk(\F_{2})=2$, and so the difference
is $rk(F_{2})-rk(H)=1$. But as we infer from Figure \ref{fig:lattice},
$\rho_{X}(H,\F_{2})=2$. Thus $\rho_{X}(H,\F_{2})>rk(\F_{2})-rk(H)$
and $H$ is \emph{not} a free factor of $\F_{2}$. As consequence,
$w$ is \emph{not} primitive in $\F_{3}$. (This example generalizes
as follows: if $w$ is a free word containing exactly $l$ different
letters, then $w$ is primitive iff we can obtain a wedge-of-loops
graph from $\G_{X}(\langle w\rangle)$ by merging $l-1$ pairs of
vertices.)
\begin{rem}
\label{rem:silva-weil2} At this point we would like to elaborate
on the differences between the algorithm presented here and the one
introduced in \cite{SW08}. Silva and Weil's presentation considers
automata and their languages. We consider the $X$-fringe $\O_{X}(H)$
and introduce the DAG $\D_{k}$ and the distance function from Definition
\ref{def:distance}. Steps 1 and 2 of our algorithm, which reduce
the problem in its very beginning to the case where $H\covers J$,
have no parallel in \cite{SW08}. Rather, they show that if $H\ff J$,
then by some sequence of {}``$i$-steps'' (their parallel of our
immediate quotients) on $H$, of length at most $rk(J)-rk(H)$, one
can obtain a core graph which is embedded in $\G_{X}(J)$ (we make
the observation that this embedded core graph can be computed in advance).
Besides shedding more light on this underlying phenomenon, our more
graph-theoretic approach has another substantial advantage: by considering
$\D_{k}$, turning the fringe $\O_{X}(H)$ into a directed graph and
stating the algorithm in the language of Theorem \ref{thr:rho=00003Drk-rk_iff_ff},
we obtain a straight-forward algorithm to identify $H$-critical subgroups
and to compute $\pi(H)$. Moreover, we obtain a straight-forward algorithm
to identify all {}``algebraic extensions'' of $H$ (Corollary \prettyref{cor:identify-alg-extensions}).
In particular, our algorithm to identify algebraic extensions substantially
improves the one suggested in \cite{KM02}, Theorem 11.3 (and see
also remark 11.4 about its efficiency). 
\end{rem}

\section{The Proof of Lemma \ref{lem:bounds_for_rho} \label{sec:The-Proof-of-bounds-lemma}}

To complete the picture, we prove the upper bound for $\rho_{X}(H,J)$
stated in Lemma \ref{lem:bounds_for_rho}. We need to show that if
$H,J\fg\F_{k}$ such that $H\covers J$, then 
\[
\rho_{X}(H,J)~~\le~~\mathrm{rk}(J)
\]

\begin{proof}
We show that $\Delta=\G_{X}(J)$ can be obtained from $\G=\G_{X}(H)$
by merging at most $\mathrm{rk}(J)$ pairs of vertices. To see this,
denote by $m$ the number of edges in $\G$, and choose some order
on these edges, $e_{1},\ldots,e_{m}$ so that for every $i$, there
is a path from the basepoint of $\G$ to $e_{i}$ traversing only
edges among $e_{1},\ldots,e_{i-1}$. (So $e_{1}$ must be incident
with the basepoint, $e_{2}$ must be incident either with the basepoint
or with the other end of $e_{1}$, etc.)

We now expose $\Delta$ step by step, each time adding the images
of the next edge of $\G$ and of its end points. Formally, denote
by $\eta$ the (surjective) morphism from $\G$ to $\Delta$, let
$\G_{i}$ be the subgraph of $\G$ that is the union of the basepoint
of $\G$ together with $e_{1},\ldots,e_{i}$ and their endpoints,
and let $\Delta_{i}=\eta(\G_{i})$. We thus have two series of subgraphs
\[
\G_{0}\subseteq\G_{1}\subseteq\ldots\subseteq\G_{m}=\G
\]
 and 
\[
\Delta_{0}\subseteq\Delta_{1}\subseteq\ldots\subseteq\Delta_{m}=\Delta
\]
 with $\Delta_{0}=\G_{0}$ being graphs with a single vertex and no
edges.

Assume that $e_{i}=(u,v)$, and w.l.o.g. that $u\in V(\G_{i-1})$.
We then distinguish between three options. A \textbf{forced} step
is when $\eta(e_{i})$ already belongs to $\Delta_{i-1}$ and then
$\Delta_{i}=\Delta_{i-1}$. A \textbf{free} step takes place when
neither $\eta(e_{i})$ nor $\eta(v)$ belong to $\Delta_{i-1}$, in
which case $\pi_{1}(\Delta_{i})=\pi_{1}(\Delta_{i-1})$. The third
option is that of a \textbf{coincidence}. This means that $\eta(e_{i})$
does not belong to $\Delta_{i-1}$ but $\eta(v)$ does. In this case,
$\Delta_{i}$ is obtained from $\Delta_{i-1}$ by connecting two vertices
by a new edge, and $\pi_{1}(\Delta_{i})$ has rank larger by $1$
from the rank of $\pi_{1}(\Delta_{i-1})$. Since the fundamental group
of $\Delta_{0}$ has rank $0$, this shows there are exactly $rk(J)$
coincidences along this process.

Assume the coincidences occurred in steps $j_{1},\ldots,j_{\mathrm{rk}\left(J\right)}$.
If $e_{j_{i}}=\left(u,v\right)$, we let $\tilde{v}\in\eta^{-1}(\eta(v))\cap V(\G_{i-1})$,
and take $\{v,\tilde{v}\}$ to be a pair of vertices of $\G$ that
we merge. (It is possible that $v=w$.) Let $w_{i}\in\F_{k}$ be be
a word corresponding to this merge in $\Gamma$. It is easy to see
by induction that $\Delta_{j_{i}}$ corresponds to the subgroup $\left\langle H,w_{1},\ldots,w_{i}\right\rangle $.
In particular, $\Delta$ corresponds to $\left\langle H,w_{1},\ldots,w_{\mathrm{rk\left(J\right)}}\right\rangle $.
We are done because all these words correspond to pairs of vertices
in $\Gamma$ (and see Remark \prettyref{remark:pairs-of-vertices}).
\end{proof}

\section{The Folding Algorithm to Construct Core Graphs \label{sec:Stallings-Folding-Algorithm}}

Finally, we present a well known algorithm to construct the core graph
of a given subgroup $H\fg\F_{k}$. The input to this process is any
finite set of words $\{h_{1},\ldots,h_{r}\}$ in the letters $\{x_{1},\ldots,x_{k}\}$
that generate $H$.

\begin{figure}[h]
\begin{centering}
\includegraphics[width=0.99\columnwidth]{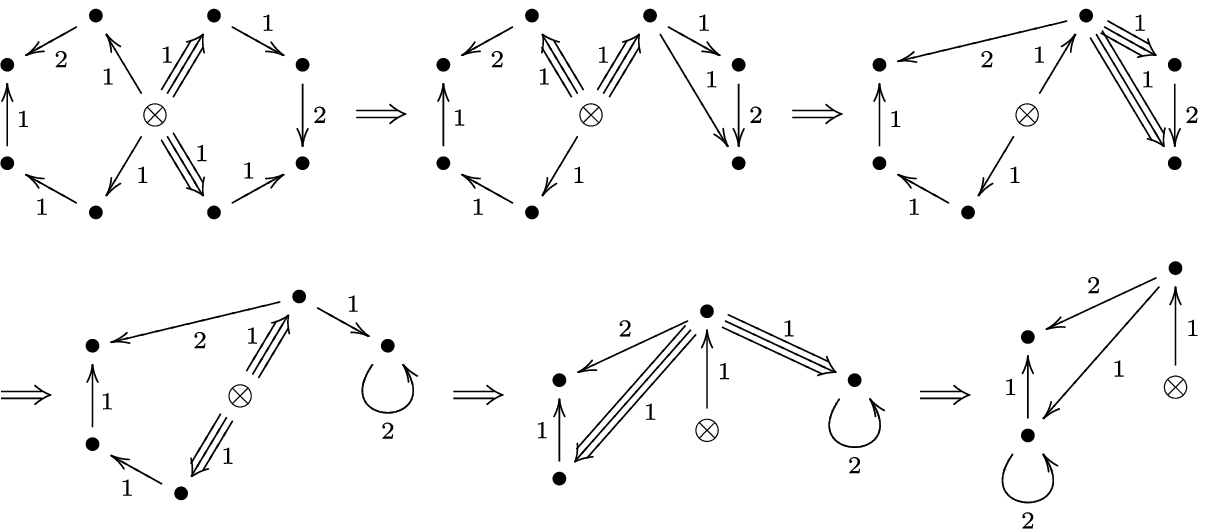}
\par\end{centering}

\caption{\label{fig:folding_process} Generating the core graph $\Gamma_{X}(H)$
of $H=\langle x_{1}x_{2}x_{1}^{-3},x_{1}^{\;2}x_{2}x_{1}^{-2}\rangle\leq\F_{2}$
from the given generating set. We start with the upper left graph
which contains a distinct loop at the basepoint for each (reduced)
element of the generating set. Then, gradually and at arbitrary order,
we merge pairs of equally-labeled edges which share the same origin
or the same terminus. One of the possible orders of merging pairs
is shown in this figure, and at each phase we mark by triple arrows
the pair of edges being merged. The graph in the bottom right is $\G_{X}(H)$,
as it has no equally-labeled edges sharing the same origin or the
same terminus.}
\end{figure}

Every element $h_{i}$ of the generating set corresponds to some path
with directed edges labeled by the $x_{i}$'s (we assume the element
is given in reduced form). Merge these $r$ paths to a single graph
by identifying all their $2r$ end-points to a single vertex which
is denoted as basepoint. Then, as long as there are two $j$-labeled
edges with the same terminus (resp. origin) for some $j$, merge the
two edges and their origins (resp. termini). Such a step is often
referred to as a \emph{Stallings' folding}. It is a fairly easy observation
that the resulting graph is indeed $\G_{X}(H)$ and that the order
of folding has no significance. To illustrate, we draw in Figure \ref{fig:folding_process}
the folding process by which we obtain the core graph $\G_{X}(H)$
of $H=\langle x_{1}x_{2}x_{1}^{-3},x_{1}^{\;2}x_{2}x_{1}^{-2}\rangle\leq\F_{2}$
from the given generating set.

\end{appendices}

\bibliographystyle{amsalpha}
\bibliography{PrimitiveBib}

\providecommand{\bysame}{\leavevmode\hbox to3em{\hrulefill}\thinspace}
\providecommand{\MR}{\relax\ifhmode\unskip\space\fi MR }
\providecommand{\MRhref}[2]{%
  \href{http://www.ams.org/mathscinet-getitem?mr=#1}{#2}
}
\providecommand{\href}[2]{#2}
\begin{thebibliography}{MVW07}

\bibitem[Abe06]{Abe06}
Mikl\'os Abert, \emph{On the probability of satisfying a word in a group},
  Journal of Group Theory \textbf{9} (2006), 685--694.

\bibitem[Bog08]{Bog08}
Oleg Bogopolski, \emph{Introduction to group theory}, EMS Textbooks in
  Mathematics, European Mathematical Society, Zurich, 2008.

\bibitem[Ger84]{Ger84}
SM~Gersten, \emph{On whitehead's algorithm}, Bull. Amer. Math. Soc., New Ser
  \textbf{10} (1984), no.~2, 281--284.

\bibitem[GS09]{GSh09}
Shelly Garion and Aner Shalev, \emph{Commutator maps, measure preservation, and
  t-systems}, Trans. Amer. Math. Soc. \textbf{361} (2009), no.~9, 4631--4651.

\bibitem[KM02]{KM02}
I.~Kapovich and A.~Myasnikov, \emph{Stallings foldings and subgroups of free
  groups}, Journal of Algebra \textbf{248} (2002), no.~2, 608--668.

\bibitem[LP10]{LP10}
Nati Linial and Doron Puder, \emph{Words maps and spectra of random graph
  lifts}, Random Structures and Algorithms \textbf{37} (2010), no.~1, 100--135.

\bibitem[LS70]{LS70}
R.C Lyndon and P.E. Schupp, \emph{Combinatorial group theory}, Springer-Verlag,
  Berlin; New York, 1970.

\bibitem[LS08]{LSh08}
Michael Larsen and Aner Shalev, \emph{Characters of symmetric groups: sharp
  bounds and applications}, Inventiones mathematicae \textbf{174} (2008),
  no.~3, 645--687.

\bibitem[LS09]{LSh09}
\bysame, \emph{Words maps and waring type problems}, J. Amer. Math. Soc.
  \textbf{22} (2009), no.~2, 437--466.

\bibitem[MVW07]{MVW07}
Alexei Miasnikov, Enric Ventura, and Pascal Weil, \emph{Algebraic extensions in
  free groups}, Geometric group theory (G.N. Arzhantseva, L.~Bartholdi,
  J.~Burillo, and E.~Ventura, eds.), Trends Math., Birkhauser, 2007,
  pp.~225--253.

\bibitem[Nic94]{Nic94}
Alexandru Nica, \emph{On the number of cycles of given length of a free word in
  several random permutations}, Random Structures and Algorithms \textbf{5}
  (1994), no.~5, 703--730.

\bibitem[PP12]{PP12}
D.~Puder and O.~Parzanchevski, \emph{Measure preserving words are primitive},
  Arxiv preprint arXiv:1202.3269 (2012).

\bibitem[RVW07]{RVW07}
A.~Roig, E.~Ventura, and P.~Weil, \emph{On the complexity of the whitehead
  minimization problem}, International journal of Algebra and Computation
  \textbf{17} (2007), no.~8, 1611--1634.

\bibitem[Seg09]{Seg09}
Dan Segal, \emph{Words: notes on verbal width in groups}, London Mathematical
  Society, Lecture note Series 361, Cambridge University Press, Cambridge,
  2009.

\bibitem[Sha09]{Sha09}
Aner Shalev, \emph{Words maps, conjugacy classes, and a non-commutative
  waring-type theorem}, Annals of Math. \textbf{170} (2009), 1383--1416.

\bibitem[Sta83]{Sta83}
John~R. Stallings, \emph{Topology of finite graphs}, Inventiones mathematicae
  \textbf{71} (1983), no.~3, 551--565.

\bibitem[Sta99]{Sta99}
\bysame, \emph{Whitehead graphs on handlebodies}, Geometric group theory down
  under (J.~Cossey, C.~F. Miller, W.D. Neumann, and M.~Shapiro, eds.), de
  Gruyter, Berlin, 1999, pp.~317--330.

\bibitem[SW08]{SW08}
P.~Silva and P.~Weil, \emph{On an algorithm to decide whether a free group is a
  free factor of another}, RAIRO - Theoretical Informatics and Applications
  \textbf{42} (2008), no.~2, 395--414.

\bibitem[Whi36a]{Whi36a}
J.H.C. Whitehead, \emph{On certain sets of elements in a free group}, Proc.
  London Math. Soc. \textbf{41} (1936), 48--56.

\bibitem[Whi36b]{Whi36b}
\bysame, \emph{On equivalent sets of elements in a free group}, Ann. of Math.
  \textbf{37} (1936), 768--800.

\bibitem[Wil98]{Wil98}
John~S. Wilson, \emph{Profinite groups}, Clarendon Press, Oxford, 1998.

\end{thebibliography}

{}
\end{document}